\documentclass[12pt]{article}
\usepackage{times}
\usepackage{natbib}
 \bibpunct{(}{)}{;}{a}{,}{,}
\usepackage{graphicx}
\usepackage{amsmath,amsfonts,amssymb,mathrsfs,amsthm,array}
\usepackage{alltt,color,url}
\usepackage[ruled,vlined]{algorithm2e}
\usepackage{xy}
\usepackage{epsfig, rotating}
\usepackage{bm}
\usepackage{rotating}
\usepackage[vlined]{algorithm2e}
\usepackage{rotating}
\usepackage{multicol}
\usepackage{titlesec}
\usepackage{latexsym}
\usepackage[pdftex,colorlinks=true,linkcolor=blue,citecolor=blue,urlcolor=blue,bookmarks=false,pdfpagemode=None]{hyperref}
\usepackage{url}
\usepackage{verbatim}
\usepackage{fancyhdr}
\usepackage{setspace}
\usepackage{paralist}
\usepackage{boxedminipage}
\usepackage{color}
\definecolor{lightgrey}{rgb}{0.85,0.85,0.85}
\definecolor{GreenYellow}   {cmyk}{0.15,0,0.69,0}
\definecolor{Yellow}        {cmyk}{0,0,1,0}
\definecolor{Goldenrod}     {cmyk}{0,0.10,0.84,0}
\definecolor{Dandelion}     {cmyk}{0,0.29,0.84,0}
\definecolor{Apricot}       {cmyk}{0,0.32,0.52,0}
\definecolor{Peach}         {cmyk}{0,0.50,0.70,0}
\definecolor{Melon}         {cmyk}{0,0.46,0.50,0}
\definecolor{YellowOrange}  {cmyk}{0,0.42,1,0}
\definecolor{Orange}        {cmyk}{0,0.61,0.87,0}
\definecolor{BurntOrange}   {cmyk}{0,0.51,1,0}
\definecolor{Bittersweet}   {cmyk}{0,0.75,1,0.24}
\definecolor{RedOrange}     {cmyk}{0,0.77,0.87,0}
\definecolor{Mahogany}      {cmyk}{0,0.85,0.87,0.35}
\definecolor{Maroon}        {cmyk}{0,0.87,0.68,0.32}
\definecolor{BrickRed}      {cmyk}{0,0.89,0.94,0.28}
\definecolor{Red}           {cmyk}{0,1,1,0}
\definecolor{OrangeRed}     {cmyk}{0,1,0.50,0}
\definecolor{RubineRed}     {cmyk}{0,1,0.13,0}
\definecolor{WildStrawberry}{cmyk}{0,0.96,0.39,0}
\definecolor{Salmon}        {cmyk}{0,0.53,0.38,0}
\definecolor{CarnationPink} {cmyk}{0,0.63,0,0}
\definecolor{Magenta}       {cmyk}{0,1,0,0}
\definecolor{VioletRed}     {cmyk}{0,0.81,0,0}
\definecolor{Rhodamine}     {cmyk}{0,0.82,0,0}
\definecolor{Mulberry}      {cmyk}{0.34,0.90,0,0.02}
\definecolor{RedViolet}     {cmyk}{0.07,0.90,0,0.34}
\definecolor{Fuchsia}       {cmyk}{0.47,0.91,0,0.08}
\definecolor{Lavender}      {cmyk}{0,0.48,0,0}
\definecolor{Thistle}       {cmyk}{0.12,0.59,0,0}
\definecolor{Orchid}        {cmyk}{0.32,0.64,0,0}
\definecolor{DarkOrchid}    {cmyk}{0.40,0.80,0.20,0}
\definecolor{Purple}        {cmyk}{0.45,0.86,0,0}
\definecolor{Plum}          {cmyk}{0.50,1,0,0}
\definecolor{Violet}        {cmyk}{0.79,0.88,0,0}
\definecolor{RoyalPurple}   {cmyk}{0.75,0.90,0,0}
\definecolor{BlueViolet}    {cmyk}{0.86,0.91,0,0.04}
\definecolor{Periwinkle}    {cmyk}{0.57,0.55,0,0}
\definecolor{CadetBlue}     {cmyk}{0.62,0.57,0.23,0}
\definecolor{CornflowerBlue}{cmyk}{0.65,0.13,0,0}
\definecolor{MidnightBlue}  {cmyk}{0.98,0.13,0,0.43}
\definecolor{NavyBlue}      {cmyk}{0.94,0.54,0,0}
\definecolor{RoyalBlue}     {cmyk}{1,0.50,0,0}
\definecolor{Blue}          {cmyk}{1,1,0,0}
\definecolor{Cerulean}      {cmyk}{0.94,0.11,0,0}
\definecolor{Cyan}          {cmyk}{1,0,0,0}
\definecolor{ProcessBlue}   {cmyk}{0.96,0,0,0}
\definecolor{SkyBlue}       {cmyk}{0.62,0,0.12,0}
\definecolor{Turquoise}     {cmyk}{0.85,0,0.20,0}
\definecolor{TealBlue}      {cmyk}{0.86,0,0.34,0.02}
\definecolor{Aquamarine}    {cmyk}{0.82,0,0.30,0}
\definecolor{BlueGreen}     {cmyk}{0.85,0,0.33,0}
\definecolor{Emerald}       {cmyk}{1,0,0.50,0}
\definecolor{JungleGreen}   {cmyk}{0.99,0,0.52,0}
\definecolor{SeaGreen}      {cmyk}{0.69,0,0.50,0}
\definecolor{Green}         {cmyk}{1,0,1,0}
\definecolor{ForestGreen}   {cmyk}{0.91,0,0.88,0.12}
\definecolor{PineGreen}     {cmyk}{0.92,0,0.59,0.25}
\definecolor{LimeGreen}     {cmyk}{0.50,0,1,0}
\definecolor{YellowGreen}   {cmyk}{0.44,0,0.74,0}
\definecolor{SpringGreen}   {cmyk}{0.26,0,0.76,0}
\definecolor{OliveGreen}    {cmyk}{0.64,0,0.95,0.40}
\definecolor{RawSienna}     {cmyk}{0,0.72,1,0.45}
\definecolor{Sepia}         {cmyk}{0,0.83,1,0.70}
\definecolor{Brown}         {cmyk}{0,0.81,1,0.60}
\definecolor{Tan}           {cmyk}{0.14,0.42,0.56,0}
\definecolor{Gray}          {cmyk}{0,0,0,0.50}
\definecolor{Black}         {cmyk}{0,0,0,1}
\definecolor{White}         {cmyk}{0,0,0,0}
\usepackage{lineno}
\usepackage[top=1in, bottom=1in, left=1in, right=1in]{geometry}

\titlespacing{\section}{0pt}{*0.9}{*0.9}
\titlespacing{\subsection}{0pt}{*0.8}{*0.8}
\titlespacing{\subsubsection}{0pt}{*0.8}{*0.8}

\newcommand{\E}{\mathcal{E}}

\def\nn{{\mathbb N}}

\def\rr{{\mathbb R}}
\newcommand{\Dim}{d}  %% Was: m
\newcommand{\Num}{p}  %% Was: d

\newcommand{\BB}[1]{{\mathbb{B}^{#1}}}
\newcommand{\tim}[1]{^{(#1)}}      % 'time', for RW's
\newcommand{\iid}{\mathrel{\mathop{\sim}\limits^{\mathrm{iid}}}}
\newcommand{\IID}{\emph{iid}} %\emph{i.i.d}
\newcommand{\Cech}{\textrm{\v{C}ech}}
\newcommand{\Alpha}{\textrm{Alpha}}
\newcommand{\Nrv}{\textrm{Nrv}}
\newcommand{\ER}{Erd{\"o}s-R{\'e}nyi}
\newcommand{\cq}[1]{[#1]}     % Clique:        [1,2,5]
\newcommand{\dist}[1]{{\textrm{dist}\left(#1\right)}}
\newcommand{\cs}[1]{\{#1\}}       % Complete set:  {1,2,5}
\newcommand{\CQ}[1]{\red{\bf\cq{#1}}}  % TRUE Clique
\newcommand{\CS}[1]{\red{\bf\cs{#1}}}  % TRUE Complete set

\newcommand{\red}[1]{\textcolor{red}{#1}}
\newcommand{\half}{{\ensuremath\textstyle\frac12}}
\newtheorem{prop}{Proposition}[section]
\newtheorem{definition}{Definition}[section]
\numberwithin{equation}{section}
\newcommand{\Ex}{\textsf{Ex}}

\newcommand{\No}{\textsf{No}}

\newcommand{\Un}{\textsf{Un}}

\newcommand\set[1]{\left\{#1\right\}} % set
\newcommand\edg[1]{\left(#1\right)}  % Graph edge-- (i,j)
\newcommand{\eg} {\textit{e.g.}}
\newcommand{\etc} {\textit{etc}}
\newcommand{\ie} {\textit{i.e.}}

\newcommand{\ooNum}{{\frac{1}{\Num}}}
\newcommand{\bV} {{\mathbf{V}}}

\newcommand{\bx} {{\mathbf{x}}}
\newcommand{\cA} {{\mathcal{A}}} % Algorithm
\newcommand{\cE} {{\mathcal{E}}}
\newcommand{\cG} {{\mathcal{G}}} % function class
\newcommand{\cH} {{\mathcal{H}}} % function class
\newcommand{\cK} {{\mathcal{K}}}
\newcommand{\cM} {{\mathcal{M}}}
\newcommand{\cP} {{\mathcal{P}}}
\newcommand{\cT} {{\mathcal{T}}}
\newcommand{\cV} {{\mathcal{V}}}

\newcommand{\cQ} {{\mathcal{Q}}}
\newcommand{\sC} {{\mathscr{C}}}
\newcommand{\sG} {{\mathscr{G}}}
\newcommand{\sL} {{\mathscr{L}}} % modified loss
\newcommand{\sP} {{\mathscr{P}}}
\newcommand{\sQ} {{\mathscr{Q}}}
\newcommand{\sS} {{\mathscr{S}}}
\def\Strut{\vrule width0pt height 10pt depth 0pt}%
\newcommand{\Alg}[1] {Algorithm~\ref{#1}}
\newcommand{\Eqn}[1] {(\ref{#1})}
\newcommand{\Fig}[1] {Figure~\ref{#1}}
\newcommand{\Figs}[2] {Figures~\ref{#1} and \ref{#2}}
\newcommand{\Sec}[1] {Section~\ref{#1}}
\newcommand{\Tab}[1] {Table~\ref{#1}}
\newcommand{\Prop}[1] {Prop.~\ref{#1}}
\newcommand{\pg}  {\textit{p}.\thinspace}

\newcommand{\hide}[1]{}
\newcolumntype{C}{>{$}c<{$}}
\newcolumntype{L}{>{$}l<{$}}
\newcolumntype{R}{>{$}r<{$}}

\begin{document}
\pagestyle{empty}

%\title{Geometric representations of distributions on hypergraphs}
\title{Geometric representations of random hypergraphs}
%\protect\thanks{Sim\'on Lunag\'omez is a Postdoctoral Fellow in the Department of Statistics at Harvard University (lunagomez@fas.harvard.edu). Sayan Mukherjee is an Associate Professor of Statistical Science, Computer Science and Mathematics at Duke University (sayan@stat.duke.edu ). Robert L.~Wolpert is a Professor of Statistical Science and the Environment at Duke University (wolpert@stat.duke.edu ). Edoardo M.~Airoldi is an Associate Professor of Statistics at Harvard University (airoldi@fas.harvard.edu). 
%%
%We are grateful to Herbert Edelsbrunner, John Harer, and Henry Wynn for helpful conversations. 
%%
%This work was partially supported 
% by the National Science Foundation under grants 
%  DMSÐ0635449, DMSÐ0732260, DMS-0757549, CAREER IIS-1149662, and IIS-1409177,
% by National Institute of Health grants 
%  NIH R01 CA123175, R01 GM096193, and NIH P50-GM081883,
% by the Army Research Office grant 
%  MURI W911NF-11-1-0036,
% and by the Office of Naval Research under grant 
%  YIP N00014-14-1-0485. 
% Edoardo M.~Airoldi is an Alfred P. Sloan Research Fellow, and a Radcliffe Institute Fellow.}}
%\author{ Sim\'on Lunag\'omez, ~Sayan Mukherjee, ~Robert L. Wolpert, ~Edoardo M.  Airoldi}
\date{}

\maketitle
\thispagestyle{empty}

%\newpage

\begin{abstract}
  We introduce a novel parametrization of distributions on hypergraphs based on the geometry of points in $\rr^\Dim$. 
  The idea is to induce distributions on hypergraphs by placing priors on point configurations via spatial processes. 
  This prior specification is then used to infer conditional independence models or Markov structure for multivariate distributions.
  %Specifically, we can recover both the junction tree factorization as well as the hyper Markov law.  
  This approach 
   supports inference of factorizations that cannot be retrieved by a graph alone, 
   leads to new Metropolis-Hastings Markov chain Monte Carlo algorithms with both local and global moves in graph space,
   and 
   generally offers greater control on the distribution of graph features than currently possible.
  We provide a comparative performance evaluation against state-of-the-art, and we
  illustrate the utility of this approach on simulated and real data.
\vfill

\end{abstract}

\singlespacing 

\newpage
\tableofcontents
% \setcounter{tocdepth}{4}
% 1 Section
% 2 Subsection
% 3 Subsubsection
% 4 Paragraph
% 5 Subparagraph

\onehalfspacing

 %%% %%% %%%
 %%% %%% %%%
 %%% %%% %%%
 
\newpage

\pagestyle{fancy}
\setcounter{page}{1}

\newpage
\section{Introduction}\label{s:intro}
% Section 1

Consider the problem of making inference on the dependence structure among random variables $X_1,...,X_p \in \mathbb{R}^p$, from $m$ replicated observations. 
% EDO : MAKE SURE THE NOTATION FOR "NUM VARIABLES" AND "NUM REPLICATES" IS CONSISTENT THROUGHOUT
The dominant formalism for this problem, is that of graphical models \citep{Laur:1996}. In this formalism. the focus is  on the first two moments of the observation vector, $X=\{X_1,...,X_p\}$, and the dependence structure is specified in terms of pairwise relations, which define an undirected graph. If such a graph is decomposable, inference is typically carried out efficiently.
Here we detail a new approach for the construction of distributions on undirected graphs, motivated by the problem of Bayesian inference of the dependence structure among random variables.  
% EDO : MOVED DOWN
%As a side
%benefit our approach also yields estimates of the conditional distributions
%given the graph.  The model space of undirected graphs grows quickly with the
%dimension of  $\set{ X_1,\ldots,X_\Num}$  (there are $2^{\Num(\Num-1)/2}$ undirected graphs on
%$\Num$ vertices) and is difficult to parametrize.  We propose a novel
%parametrization and a simple, flexible family of prior distributions on $\cG$
%and on Markov probability distributions with respect to $\cG$ \citep
%{Dawi:Laur:1993}; this parametrization is based on computing the intersection
%pattern of a system of convex sets in $\rr^\Dim$.
%% and that is simple and flexible enough to place meaningful informative
%% prior distributions on this model space.support meaningful expressions
%% Inference of the factorization of a multivariate distribution can be
%% separated into two parts: inference of the graph $\cG$ and inference
%% of the marginal distributions given by the graph using prior
%% specification
%The novelty and main contribution of this paper is structural inference for
%graphical models, specifically, the proposed representation of graph spaces
%allows for flexible prior distributions and new Markov chain Monte Carlo
%(MCMC) algorithms.

\subsection{Related work}

It is common to model the joint probability distribution of a family of
$p$ random variables $\{ X_1,\ldots,X_p\}$ in two stages. First 
specify the \emph{conditional dependence structure} of the distribution, then
 specify details of the {conditional distributions} of the variables within
that structure \citetext{see \pg1274 of \citealt{Dawi:Laur:1993}, or \pg180
  of \citealt {Besa:1975}, for example}.  The structure may be summarized in
a variety of ways in the form of a graph $\cG=(\cV,\cE)$ whose vertices
$\cV=\{1,...,p\}$ index the variables $\set{X_i}$ and whose edges
$\cE\subseteq \cV\times\cV$ in some way encode conditional dependence.
%
%In the common DAG or Directed Acyclic Graph approach each edge $(i,j) \in\cE$
%suggests some direct influence of $X_i$ on $X_j$, and the full distribution
%is then written as a product of
%%independent marginal distributions for the ``orphan'' variables $\cO=
%%\{j:~\nexists (i,j)\in\cE\}$ and conditional distributions for each
%%non-orphan $j\in\cV\setminus\cO$, 
%conditional distributions for each $X_j$ given its ``parents'' $p(j)=\{i\in
%\cV:~(i,j)\in\cE\}$.  We follow a different approach, suitable for joint
%distributions with an everywhere-positive density function
We follow the Hammersley-Clifford approach \citep {Besa:1974,Hamm:Clif:1971},
in which $(i,j)\in\cE$ if and only if the conditional distribution of $X_i$
given all other variables $\set{X_k:~k\ne i}$ depends on $X_j$, \ie, differs
from the conditional distribution of $X_i$ given $\set{X_k:~k\ne i,j}$.  In
this case the distribution is said to be Markov with respect to the graph.
One can show that this graph is symmetric or \emph{undirected}, \ie, all the
elements of $\cE$ are unordered pairs.
% (so $(i,j)\in\cE\Leftrightarrow(j,i)\in\cE$). 
%$(i,j)\in\cE$ if and only if $(j,i)\in\cE$.  

The simultaneous inference of a decomposable graph and marginal distributions
in a fully Bayesian framework was approached in \citep {Gree:1995} using
local proposals to sample graph space.  A promising extension of this
approach called Shotgun Stochastic Search (SSS) takes advantage of parallel
computing to select from a batch of local moves \citep {Jone:Carv:etal:2005}.
A stochastic search method that incorporates both local moves and more
aggressive global moves in graph space has been developed by \citet
{Scot:Carv:2008}.  These stochastic search methods are intended to identify
regions with high posterior probability, but their convergence properties are
still not well understood.  Bayesian models for non-decomposable graphs have
been proposed by \citet {Rove:2002} and by \citet*{Wong:Cart:Kohn:2003}.
These two approaches focus on Monte Carlo sampling of the posterior
distribution from specified hyper Markov prior laws.  Their emphasis is on
the computational problem of Monte Carlo simulation, not on that of
constructing interesting informative priors on graphs.  We think there is
need for methodology that offers both efficient exploration of the model
space and a simple and flexible family of distributions on graphs that can
reflect meaningful prior information.

\ER\ random graphs (those in which each of the $\binom{p}{2}$ possible
undirected edges $(i,j)$ is included in $\cE$ independently with some
specified probability $\alpha \in[0,1]$), and variations where the edge inclusion
probabilities $\alpha_{ij}$ are allowed to be edge-specific, have been used to
place informative priors on decomposable graphs \citep {Heck:Geig:Chic:1995,
  Mans:Kemp:etal:2006}.  The number of parameters in this prior specification
can be enormous if the inclusion probabilities are allowed to vary, and some
interesting features of graphs (such as decomposability) cannot be expressed
solely through edge probabilities.  \citet {Mukh:Spee:2008} developed methods
for placing informative distributions on directed graphs by using
\emph{concordance functions} (functions that increase as the graph agrees
more with a specified feature) as potentials in a Markov model.
% This is a conceptually sound approach, however 
This approach is tractable, but it is still not clear how to encode certain
common assumptions within such a framework.

For the special case of jointly Gaussian variables $\set{X_j}$, or those with
arbitrary marginal distributions $F_j(\cdot)$ whose dependence is adequately
represented in Gaussian copula form $X_j=F_j^{-1}\big (\Phi(Z_j)\big)$ for
jointly Gaussian $\set{Z_j}$ with zero mean and unit-diagonal covariance
matrix $C$, the problem of studying conditional independence reduces to
a search for zeros in the precision matrix $C^{-1}$.  This approach \citep
[see][for example] {Hoff:2007} is faster and easier to implement than ours in
cases where both are applicable, but is far more limited in the range of
dependencies it allows.  For example, a three-dimensional model in which each
pair of variables is conditionally independent given the third cannot be
distinguished from a model with complete joint dependence of the three
variables (we return to this example in \Sec{Sec:FacNerve}).

\subsection{Contributions}

In this article we establish a novel approach to parametrize spaces of
graphs.  For any integers $p,\Dim\in\nn$, we show in \Sec{s:rgg} how to
use the geometrical configuration of a set $\set{v_i}$ of $p$ points in
Euclidean space $\rr^\Dim$ to determine a graph $\cG=(\cV,\cE)$ on
$\cV=\set{v_1,...,v_p}$.  Any prior distribution on point sets $\set{v_i}$
induces a prior distribution on graphs, and sampling from the posterior
distribution of graphs is reduced to sampling from spatial configurations of
point sets--- a standard problem in spatial modeling.  Relations between
graphs and finite sets of points have arisen earlier in the fields of
computational topology \citep {Edel:Hare:2008} and random geometric graphs
\citep {Penr:2003}.  From the former we borrow the idea of \emph{nerves},
\ie, simplicial complexes computed from intersection patterns of convex
subsets of $\rr^\Dim$; the $1$-skeletons (collection of $1$-dimensional
simplices) of nerves are geometric graphs.

As a side
benefit our approach also yields estimates of the conditional distributions
given the graph.  The model space of undirected graphs grows quickly with the
dimension of  $\{ X_1,\ldots,X_p\}$  (there are $2^{p(p-1)/2}$ undirected graphs on
$p$ vertices) and is difficult to parametrize.  We propose a novel
parametrization and a simple, flexible family of prior distributions on $\cG$
and on Markov probability distributions with respect to $\cG$ \citep
{Dawi:Laur:1993}; this parametrization is based on computing the intersection
pattern of a system of convex sets in $\rr^\Dim$.
% and that is simple and flexible enough to place meaningful informative
% prior distributions on this model space.support meaningful expressions
% Inference of the factorization of a multivariate distribution can be
% separated into two parts: inference of the graph $\cG$ and inference
% of the marginal distributions given by the graph using prior
% specification
The novelty and main contribution of this paper is structural inference for
graphical models, specifically, the proposed representation of graph spaces
allows for flexible prior distributions and new Markov chain Monte Carlo
(MCMC) algorithms.

From the random geometric graph approach we gain understanding about the
induced distribution on graph features when making certain features of a
geometric graph (or hypergraph) stochastic.

\section{Background and preliminaries}

\subsection{Graphical models}\label{sec:graphical models}

The graphical models framework is concerned with the representation of
conditional dependencies for a multivariate distribution in the form of a
graph or hypergraph.  We first review relevant graph theoretical concepts and
then relate these concepts to factorizing distributions.

A \emph{graph} $\cG$ is an ordered pair $(\cV,\cE)$ of a set $\cV$ of
\emph{vertices} and a set $\cE\subseteq\cV\times\cV$ of \emph{edges}.  If all
edges are unordered (resp., ordered), the graph is said to be
\emph{undirected} (resp., \emph{directed}).  All graphs considered in this
paper are undirected, unless stated otherwise.  A \emph {hypergraph}, denoted
$\cH$, consists of a vertex set $\cV$ and a collection $\cK$ of unordered
subsets of $\cV$ (known as \emph{hyperedges}); a graph is the special case
where all the subsets are vertex pairs.  A graph is \emph {complete} if
$\cE=\cV\times\cV$ contains all possible edges; otherwise it is \emph
{incomplete}.  A complete subgraph that is maximal with respect to inclusion
is a \emph {clique}.  Denote by $\sC(\cG)$ and $\sQ(\cG)$, respectively, the
collection of complete sets and cliques of $\cG$.  A \emph{path} between two
vertices $\set{v_i,v_j}\in \cV$ is a sequence of edges
% and vertices \uhoh{why vertices?} 
connecting $v_i$ to $v_j$.  A graph such that any pair of vertices can be
joined by a unique path is a \emph{tree}.  A \emph{decomposition} of an
incomplete graph $\cG=(\cV,\cE)$ is a partition of $\cV$ into disjoint
nonempty sets $(A,B,S)$ such that $S$ is complete in $\cG$ and
\textit{separates} $A$ and $B$, \ie, any path from a vertex in $A$ to a
vertex in $B$ must pass through $S$.  Iterative decomposition of a graph
$\cG$ such that at each step the separator $S_i$ is minimal and the subsets
$A_i$ and $B_i$ are nonempty generates the \emph{prime components} of $\cG$,
the collection of subgraphs that cannot be further decomposed.  If all
prime components of a graph $\cG$ are complete, then $\cG$ is said to be
\emph{decomposable}.  Any graph $\cG$ can be represented as a tree
$\cT$ whose vertices are its prime components $\sP(\cG)$; this is called its
\emph{junction tree} representation.  A junction tree is a hypergraph.

Let $\cP$ be a %n absolutely continuous
probability distribution on $\rr^p$
% (\ie, one with a density function) 
and $X = (X_1,\ldots,X_p)$ a random vector with distribution $\cP$.
\emph{Graphical modeling} is the representation of the Markov or conditional
dependence structure among the components $\set{X_i}$ in the form of a graph
$\cG=(\cV,\cE)$.  Denote by $f(x)$ the joint density function of $\set{X_i}$
(or probability mass function for discrete distributions--- more generally,
density for an arbitrary reference measure).  The distribution $\cP$ (and
hence its density $f(x)$) may depend implicitly on a vector $\theta$ of
parameters, taking values in some set $\Theta_\cG$, which in some cases will
depend on the graph $\cG$; write $\Theta=\sqcup\Theta_\cG$ for the disjoint
union of the parameter spaces for all graphs on $\cV$.

Each vertex $v_i\in \cV$ is associated with a variable $X_i$, and the edges
$\cE$ determine how the distribution factors.  The density $f(x)$ for the
distribution can be factored in a variety of ways associated with the graph
$\cG$ \citep[\pg35] {Laur:1996}.  It may be factored in terms of complete
sets $a \in \sC(\cG)$:
\begin{subequations}\label{eq:factors}
\begin{align}
  f(x)&=\prod_{a \in \sC(\cG)} \phi_a(x_a\mid \theta_a),\label{eq:fac-comp}\\
\intertext{or similarly in terms of cliques $a\in\sQ$ (assuming $f$ is positive, 
according to the Hammersley-Clifford theorem); if $\cG$ is
  decomposable then $f(x)$ may also be factored in junction-tree form as:}
% f(x)&=\prod_{a \in \sQ(\cG)} \psi_a(x_a\mid \theta_a),\label{eq:fac-cliq}\\
  f(x)&=\frac{\prod_{a \in \sP(\cG)} \psi_a(x_a\mid \theta_a)}
        {\prod_{b \in \sS(\cG)} \psi_b(x_b\mid \theta_b)}, \label{eq:fac-jt}
\end{align}
\end{subequations}
where $\sP(\cG)$ and $\sS(\cG)$ denote the prime factors and separators of
$\cG$, respectively, and where $\psi_a(x_a\mid \theta_a)$
%(resp, $\psi_b(x_b\mid \theta_b)$) denotes the marginal joint densities of
%the components $x_a$ (resp, $x_b$) for prime factors $a\in\sP(\cG)$ (resp,
%separators $b\in\sS(\cG)$) 
denotes the marginal joint density for the components $x_a$ for prime factors
$a\in\sP(\cG)$ and $\psi_b(x_b\mid\theta_b)$ that for separators
$b\in\sS(\cG)$ \citep [Eqn.~(6)] {Dawi:Laur:1993}.  In the Gaussian case, a
similar factorization to \Eqn {eq:fac-jt} holds even for non-decomposable
graphs \citep[Prop.~2] {Rove:2002}.

The prior distributions required for Bayesian inference about models of the
form \eqref {eq:factors} may be specified by giving a marginal distribution
on the set of all graphs $\cG \in \sG_\Num$ on $\Num$ vertices and conditional
distributions on each $\Theta_\cG$, the space of parameters for that graph:
\begin{equation}\label{eq:ThreeStages}
%p(X\mid \theta,\cG), \quad   % This is "f"
  p(\cG, \theta) = p(\cG) ~ p(\theta \mid \cG),  %\textrm{and} \quad
  \qquad \cG\in\sG_\Num,~\theta\in\Theta_\cG
\end{equation}
where $\theta\in\Theta_\cG$ determines the parameters $\set{\theta_a : a\in
\sC(\cG)}$ or $\set{\theta_a : a\in \sP(\cG)}$ and $\set{\theta_b : b\in
\sS(\cG)}$.  \citet {Giud:Gree:1999} pursue this approach in the Gaussian
case, while \citet {Dawi:Laur:1993} offer a rigorous framework for specifying
more general prior distributions on $\Theta_\cG$.
% the parameters involved in the factors given the graph.
Such priors, called \emph {hyper Markov laws}, inherit the conditional
independence structure from the sampling distribution, now at the parameter
level.  The hyper Inverse Wishart, useful when the factors are multivariate
normal, is by far the most studied hyper Markov law.  Most previously studied
models of the form \eqref {eq:ThreeStages} specify very little structure on
$p(\cG)$ \citep{Giud:Gree:1999,Heck:Geig:Chic:1995, Rove:2002}--- typically
$p(\cG)$ is taken to be a uniform distribution on the space of decomposable
(or unrestricted) graphs, or perhaps an \ER\ prior to encourage sparsity
\citep{Mans:Kemp:etal:2006}, with no additional structure or constraints and
hence no opportunity to express prior knowledge or belief.
 
Two inference problems arise for the model specified in \eqref
{eq:ThreeStages}: inference of the entire joint posterior distribution of the
graph and factor parameters, $(\theta,\cG)$, or inference of only the \emph
{conditional independence structure}, which entails comparing different
graphs via the marginal likelihood
\begin{equation}\label{eq:marg likeli}\notag
\Pr\set{\cG\mid x}\propto\int_{\Theta_\cG}\,
   f(x\mid \theta,\cG)\,p(\cG)\, p(\theta \mid \cG) \,d\theta.
\end{equation}
%_{\mbox{\tiny Supp}(\cG)}
%where $\mbox{Supp}(\cG)$ denotes the support of $p(\theta \mid \cG)$.
Inference about $\cG$ may now be viewed as a Bayesian model selection
procedure \citep [see][\pg348] {Robe:2001}.

\subsection{Geometric graphs}\label {s:rgg}
% Section 2

Most methodology for structural inference in graphical models either assumes
little prior structure on graph space, or else represents graphs using high
dimensional discrete spaces with no obvious geometry or metric.  In either
case prior elicitation and posterior sampling can be challenging.  In this
section we propose parametrizations of graph space that will be used in
\Sec{Sec:RandomGeometric} to specify flexible prior distributions and to
construct new Metropolis\slash Hastings MCMC algorithms
%\citep [Ch.~7] {Hast:1970, Robe:Case:2004} 
with local and global moves.  The key idea for this parametrization is to
construct graphs and hypergraphs from intersections of convex sets in
$\rr^\Dim$.

\begin{figure}[hbt]
\begin{center}
\includegraphics[height=5in]{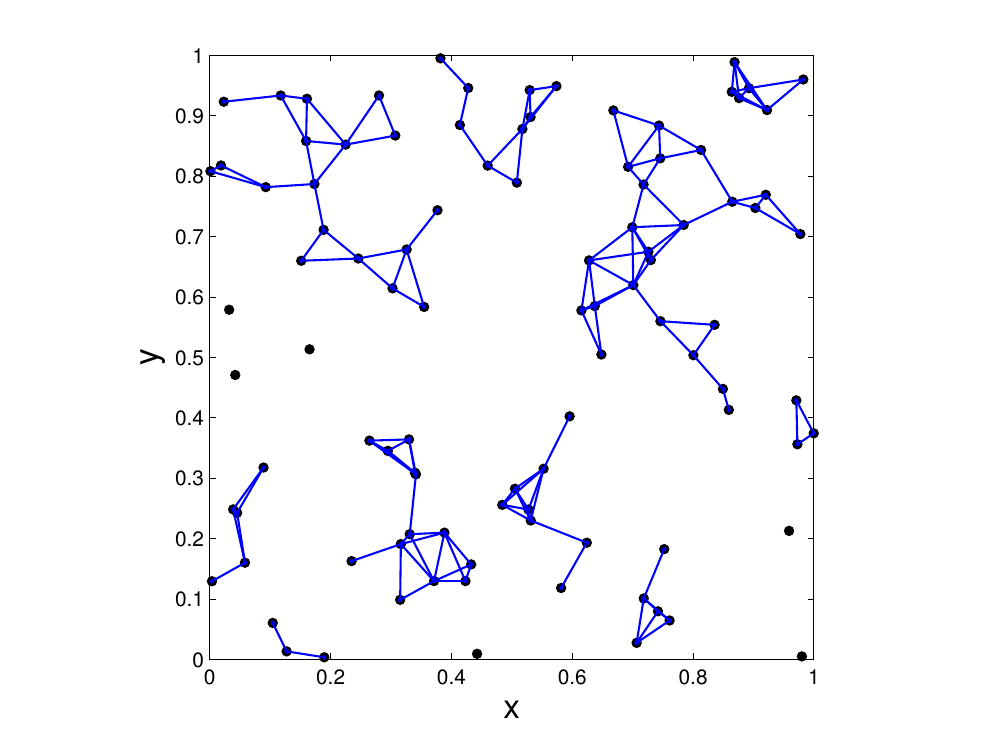}
\caption{Proximity graph for $100$ vertices and radius
  $r=0.05$.\label{fig:RGG}} 
\end{center}
\end{figure}

\begin{figure}[hbt]
\begin{center}
\includegraphics[height=120mm]{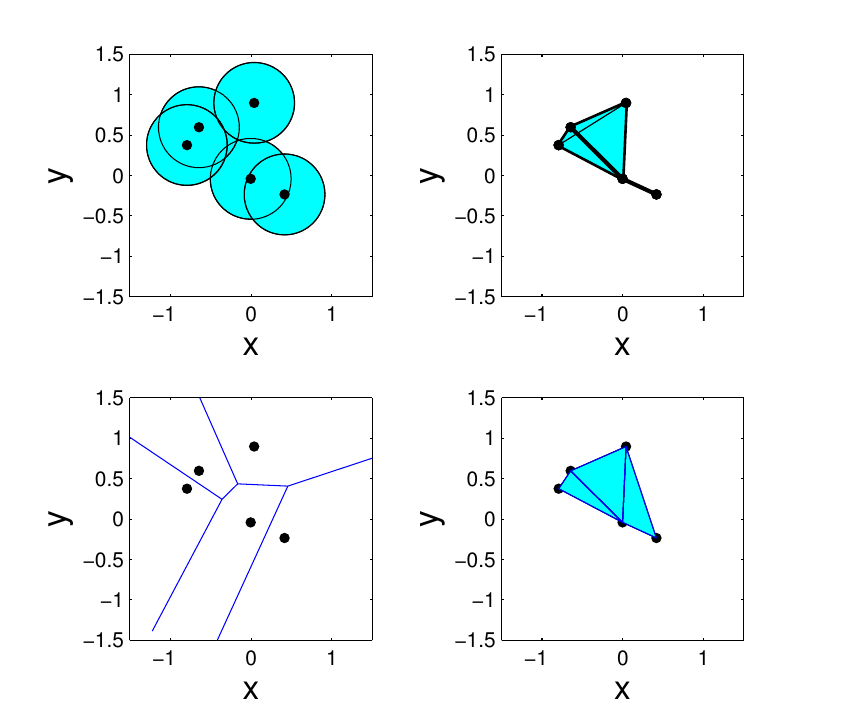}
\caption{(a) A set of vertices in $\rr^2$ are used to construct a
  family of disks of radius $r=0.5$.  (b) The nerve of this convex
  set.  This is an example of a \Cech\ complex.  (c) For the same
  vertex set the Voronoi diagram is computed.  (d) The nerve of
  the Voronoi cells is obtained.  This is an example of the
  Delaunay triangulation.  Note that the maximum clique size of the Delaunay
  is bounded by the dimension of the space of the vertex set plus one; such a
  restriction does not apply to the \Cech\ complex.}
\label{fig:Cech_Delaunay}
\end{center}
\end{figure}

We illustrate the approach with an example.  Fix a convex region $A\subset
\rr^\Dim$ and let $\cV\subset A$ be a finite set of $\Num$ points.  For each number
$r\ge0$, the proximity graph $\mbox{Prox} (\cV,r)$ (see \Fig {fig:RGG}) is
formed by joining every pair of (unordered) elements in $\cV$ whose distance
is $2r$ or less, \ie, whose closed balls of radius $r$ intersect.  As $r$
ranges from $0$ to half the diameter of $A$, the graph $\mbox{Prox}(\cV,r)$
ranges from the totally disconnected graph to the complete graph.
%  This graph, parameterized by the set of points $\cV$ and radius $r$, and
%  the geometry of the point set is used explicitly.  
This example is a particular case of a more general construction illustrated
in \Fig {fig:Cech_Delaunay}; hypergraphs can be computed from properties of
intersections of classes of convex subsets in Euclidean space.  The convex
sets we consider are subsets of $\rr^\Dim$ that are simple to parametrize and
compute.  The key concept in our construction is the \emph{nerve}:
\begin{definition}[Nerve]\label{def:nerve}
  Let $F=\set{A_j,\ j \in I}$ be a finite collection of distinct nonempty
  convex sets.  The \emph{nerve} of $F$ is given by
\[
    \Nrv(F)=\set{ \sigma \subseteq I :
                   \bigcap_{j \in \sigma}A_j\neq \varnothing }.
\]
\end{definition}
\noindent
The nerve of a family of sets uniquely determines a hypergraph.  We use the
following three nerves in this paper to construct hypergraphs \citep[for more
details, see][]{Edel:Hare:2008}.
%\uhoh{A nerve ``is'' a hypergraph or ``used to construct'' a hypergraph?}
%These constructions are used extensively in computational topology and
%geometry \citep{Edel:Hare:2008}.

\begin{definition}[\Cech\ Complex]
  Let $\cV$ be a finite set of points in $\rr^\Dim$ and $r>0$.  Denote by
  $\BB \Dim$ the closed unit ball in $\rr^\Dim$.  The \emph{\Cech\ complex}
  corresponding to $\cV$ and $r$ is the nerve of the sets $B_{v,r}=v+r\BB
  \Dim$, $v\in \cV$.  This is denoted by $\Nrv(\cV,r,\Cech)$.
\end{definition}

\begin{definition}[Delaunay Triangulation]
  Let $\cV$ be a finite set of points in $\rr^\Dim$.  The \emph{Delaunay
    triangulation} corresponding to $\cV$ is the nerve of the sets
  $C_{v}=\set{x\in \rr^\Dim:~ \|x-v\|\leq \|x-u \|, \ u\in \cV}$ for $v\in \cV$.
  This is denoted by $\Nrv(\cV,\text{Delaunay})$, and the sets $C_v$ are
  called \emph{Voronoi cells}.
\end{definition}

\begin{definition}[Alpha Complex]
  Let $\cV$ be a finite set of points in $\rr^\Dim$ and $r>0$.  The \emph {Alpha
    complex} corresponding to $\cV$ and $r$ is the nerve of the sets
  $B_{v,r}\cap C_{v}$, $v\in \cV$.  This is denoted by $\Nrv(\cV,r,\Alpha)$.
\end{definition}

\begin{figure}[t]
\begin{center}
  \includegraphics[height=80mm]{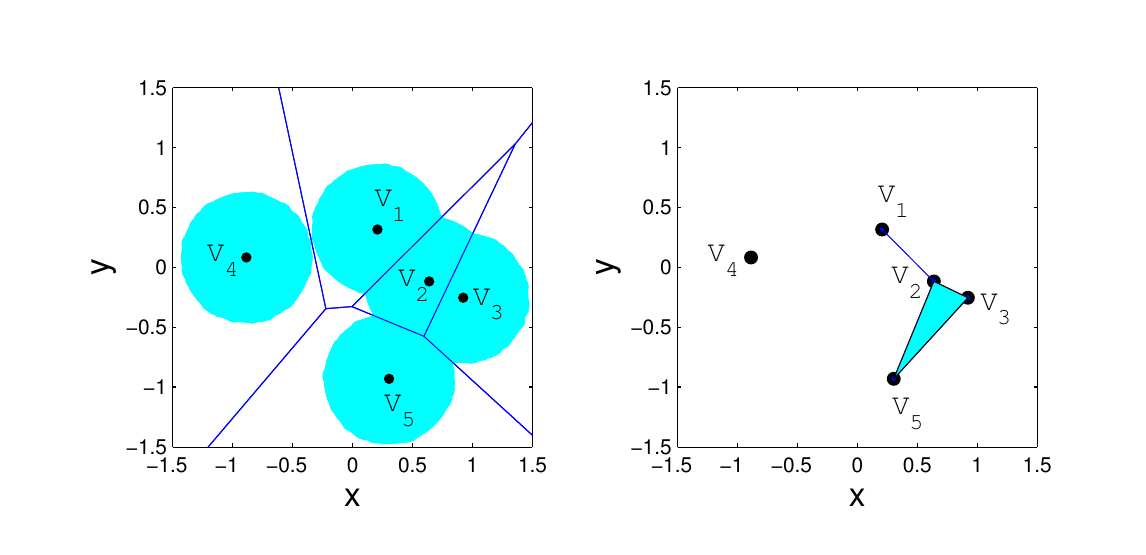}
  \caption{(a) Intersection of balls and Voronoi cells computed using $r=0.5$ and the vertex set listed in \Tab {tab:VertexExample}. (b) The corresponding Alpha complex.}\label{fig:Alpha_example}
\end{center}
\end{figure}

Here we illustrate the idea of nerve and specifically, the idea of alpha complex. Consider the vertex set displayed in Table  \ref{tab:VertexExample} and $r=0.5$. The sets $B_{v,r}\cap C_{v}$ and the corresponding nerve (alpha complex) are illustrated in Figure \ref{fig:Alpha_example}. Since the set indexed by $V_4$ does not intersect with any other $B_{v,r}\cap C_{v}$, it will produce an isolated vertex in the nerve. The set indexed by $V_1$ only intersects with the set indexed by $V_2$, therefore there will be an edge joining $V_1$ and $V_2$ in the nerve. $V_2$, $V_3$ and $V_5$ intersect as pairs, therefore, the edges of the triangle with vertices $V_2$, $V_3$ and $V_5$ will be in the nerve. Since the sets indexed by $V_2$, $V_3$ and $V_5$ also intersect as a triad, the facet or face of the triangle is also included in the nerve.

\begin{table}[t]
\begin{center}
\begin{tabular}{|c|RRRRR|}
\hline\Strut
Coordinate & V_1\quad{}&  V_2\quad{}&  V_3\quad{}&  V_4\quad{} & V_5\quad{} \\
\hline\Strut
 $x$ & 0.2065 &   0.6383 &   0.9225 &  -0.8863 &  0.3043 \\
 $y$ & 0.3149 &  -0.1193 &  -0.2544 &   0.0816 & -0.9310 \\
\hline
\multicolumn{6}{c}{}
\end{tabular}
\caption{Vertex set used for generating a family of sets and the corresponding nerve.}\label{tab:VertexExample}
\end{center}
\end{table}

The nerve of a family of sets is a particular class of hypergraphs known as
(abstract) simplicial complexes.
%\uhoh{old def was wrong--- $\cK$ isn't a subset of $\cV$, it's a subset of
%$2^{\cV}$}
\begin{definition}[Abstract simplicial complex]
  Let $\cV$ be a finite set.  A \emph{simplicial complex} with base set $\cV$
  is a family $\cK$ of subsets of $\cV$ such that $\tau\in\cK$ and $\sigma
  \subseteq \tau$ %for $\sigma,\tau\subset\cV$ 
  implies $\sigma \in \cK$.  The elements of $\cK$ are called simplices, and
  the number of connected components of $\cK$ is denoted $\sharp(\cK)$.
\end{definition}

%\uhoh{Do we really need \emph{abstract} simplicial complexes? Also, seems
%  misleading to use $\cV$ for the ``abstract set'' here and for vertices
%  elsewhere.  Usual def is something like:}
%\begin{definition}[Simplicial Complex]
%  A \emph{simplex} is the convex hull of a set of $(n{+}1)$ affinely
%  independent points in Euclidean space of dimension $n$ or greater.  The
%  $(n+1)$ points are the \emph{vertices} of the simplex, and the convex hull
%  of any subset of the vertices of a simplex is a \emph {face}.  A
%  \emph{simplicial complex} is a set $\cK$ of simplices with the property
%  that any face of a simplex $\sigma\in\cK$ is also in $\cK$, and that the
%  intersection $\sigma_1\cap\sigma_2$ of any two simplices $\sigma_1,
%  \sigma_2\in\cK$ is a face of both $\sigma_1$ and $\sigma_2$.  The
%  \emph{rank} of a simplex $\sigma\in\cK$ is the number of distinct vertices
%  in $\sigma$, or one plus the dimension; this is the maximal length of a
%  chain $\varnothing=\sigma_0\subset\sigma_1\subset \ldots
%  \subset\sigma_{|\sigma|}=\sigma$ with each $\sigma_j\in\cK$.
%\end{definition}

The nerve of a collection of sets is always a hypergraph; in simple
cases, only vertex pairs arise so the $1$-skeleton determines a unique graph.

\begin{definition}[$p$-skeleton]
  Let $\cK$ be a simplicial complex, and denote by $|\tau|$ the cardinality
  of a simplex $\tau\in\cK$.  The \emph{$p$-skeleton} of $\cK$ is the
  collection of all $\tau \in \cK$ such that $|\tau|\leq p+1$.  The elements
  of the $p$-skeleton are called $p$-simplices and the $1$-skeleton is just a
  graph (more precisely, it is $\cV\cup\cE$ for a uniquely determined graph
  $\cG=(\cV,\cE)$).
\end{definition}

\begin{figure}[b!]
\begin{center}
\includegraphics[width=\textwidth]{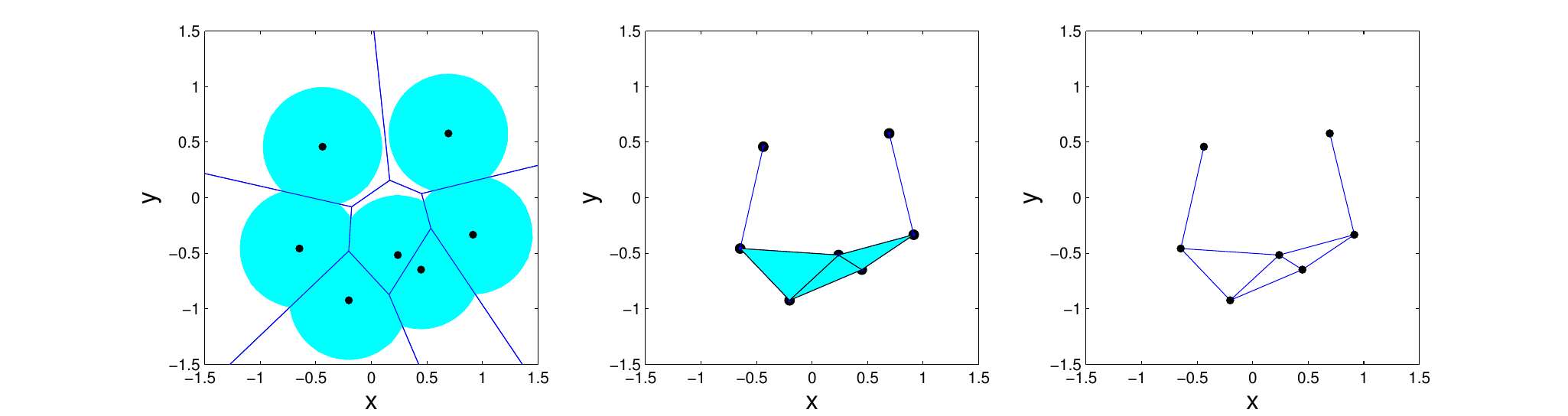}
\caption{(a) Given a set of vertices and a radius ($r=0.5$) one can compute $A_i=C_i\cap B_i$, where $C_i$ is the Voronoi cell for vertex $i$ and $B_i$ is the ball of radius $r$ centered at vertex $i$. (b) The Alpha complex is the nerve of the $A_i$'s. (c) Often the main interest will be the $1$-skeleton of the complex, which is the subset of the nerve that corresponds to (nonempty) pairwise intersections.}\label{fig:Alpha_Complex}
\end{center}
\end{figure}

The $1$-skeleton of a nerve is the graph obtained by considering only
nonempty \emph{pairwise} intersections.  The process of obtaining the nerve
and the $1$-skeleton from a family of sets is illustrated in
\Fig{fig:Alpha_Complex}.  Different families of convex sets in $\rr^\Dim$
induce different restrictions in graph space: for the Delaunay triangulation
and the Alpha complex, for example, clique sizes cannot exceed $\Dim+1$.
Although no such blanket restriction applies to the \Cech\ complex, for this
complex some graphs are still unattainable--- for example, no \Cech\ complex
can include a star graph whose central node has degree higher than the
``kissing number,'' \ie, maximal number of disjoint unit hyperspheres
touching a given hypersphere, 6 for $\Dim=2$, 12 for $\Dim=3$, \etc.
% $\Dim=2$ dimensional \Cech\ complex can include a star graph of degree
% exceeding six.

The \Cech\ and Alpha complexes are hypergraphs indexed by a finite set
$\cV=\set{V_1,\ldots,V_\Num }\subset\rr^\Dim$ and a size parameter $r\ge0$.  Each
induces a parametrization on the space of hypergraphs $(\cV,r)\mapsto
\cH(\cV,r)$.  The class $\cA$ of convex sets used to compute the nerve
determines the space of hypergraphs.  To keep the notation simple, $\cA$ will
be implicit whenever obvious by the context.  We will use $A(\cV,r)$ to denote
a generic element of $\cA$ for either the \Cech\ or the Alpha complex.
Similarly, $1$-skeletons of nerves induce a parametrization of the spaces of
graphs $(\cV,r)\mapsto \cG(\cV,r)$.

Two principal advantages of this approach are:
\begin{enumerate}
\item For each family of convex sets $\set{\cA}$, the number of parameters
  needed to specify the graph $\cG$ or hypergraph $\cH$ grows only linearly
  with the number of vertices;
\item The hypergraph parameter space will be a subset of $\rr^\Dim$, a very
  convenient parameter space for MCMC sampling.
\end{enumerate}

\subsection{Random geometric graphs}
\label{Sec:RandomGeometric}

% Section 3

In \Sec{s:rgg} we demonstrated how the geometry of a set $\cV$
of $\Num$ points in $\rr^\Dim$ can be used to induce a graph $\cG$.  In this
section we explore the relation between prior distributions on \emph{random}
sets $\bV$ of points in $\rr^\Dim$ and features of the induced distribution on
graphs $\cG$, with the goal of learning how to tailor a point process model
to obtain graph distributions with desired features.
%  The following construction will be used to study these questions.
 
\begin{definition}[Random Geometric Graph]\label{Def:RGG}
  Fix integers $\Num,\Dim\in\nn$ and let $\bV=(V_1,\ldots,V_\Num)$ be drawn from a
  probability distribution $\cQ$ on $(\rr^\Dim)^\Num$.  For any class $\cA$ of
  convex sets in $\rr^\Dim$ and radius $r>0$, the graph $\cG(\bV,r,\cA)$ is said
  to be a \emph{Random Geometric Graph} (RGG).
\end{definition}

While Definition \ref{Def:RGG} is more general than that of \cite[\pg2]
{Penr:2003}, it still cannot describe all the random graphs discussed in
\citep {Penr:Yuki:2001} (for example, those based on $k$-neighbors cannot in
general be generated by nerves).  For $\cA$ we will use closed balls in
$\rr^\Dim$ or intersections of balls and Voronoi cells; most often $\cQ$ will be
a product measure under which the $\set{V_i}$ will be $\Num$ independent
identically distributed (\IID) draws from some marginal distribution $\cQ_M$
on $\rr^\Dim$, such as the uniform distribution on the unit cube $[0,1]^\Dim$ or
unit ball $\BB \Dim$,
%Typically $\bV=(V_1,\ldots, V_\Num)$ will have a distribution $\cQ$ under which
%$\set{V_i} \iid \cQ_M$ for some marginal distribution $\cQ_M$ on $\rr^\Dim$, 
but we will also explore the use of repulsive processes for $\bV$ under which
the points $\set{V_i}$ are more widely dispersed than under independence.
It is clear that different choices for $\cA$, $\cQ$ and $r$ will have an
impact on the support of the induced RGG distribution.  To make this notion
precise we define feasible graphs. 
 
\begin{definition}[Isomorphic]
  Write $\cG_1\cong\cG_2$ for two graphs $\cG_i= (\cV_i,\cE_i)$ and call the
  graphs \emph{isomorphic} if there is a 1:1 mapping $\chi:\cV_1\to\cV_2$
  such that $\edg{v_i,v_j}\in\cE_1 \Leftrightarrow \edg{\Strut
    \chi(v_i),\chi(v_j)} \in\cE_2$ for all $v_i,v_j\in\cV_1$.
\end{definition}

\begin{definition}[Feasible Graph]\label{d:feas}
  Fix numbers $\Dim,\Num\in\nn$, a class $\cA$ of convex sets in $\rr^\Dim$, and a
  distribution $\cQ$ on the random vectors $\bV$ in
  $(\rr^\Dim)^\Num$.  A graph $\Gamma$ is said to be \emph{feasible} if for some
  number $r>0$,
\begin{displaymath}
  \Pr\set{\cG(\bV,r,\cA) \cong \Gamma}>0.
\end{displaymath}
\end{definition}

In contrast to \ER\ models, where the inclusion of graph edges are
independent events, the RGG models exhibit edge dependence that depends on
the metric structure of $\rr^\Dim$ and the class $\cA$ of convex sets used to
construct the nerves.
%  The metric structure of the space can ensure that if vertex $x$ is
%  adjacent to $y$ and $y$ is adjacent to $z$ then $x$ and $z$ will be
%  adjacent.  \uhoh{How?  That sounds wrong, by triangle inequality.  The
%  structure can make that \emph{likely}... is that what was meant?  Or am I
%  missing something?}  Using the Delaunay triangulation to compute the
%  nerves results in the existence of an edge between two vertices depending
%  on the position of all the vertices and not just the two vertices in
%  question.
 
There is an extensive literature describing asymptotic distributions for a
variety of graph features such as: subgraph counts, vertex degree, order of
the largest clique, and maximum vertex degree \cite[for an encyclopedic
account of results for the important case of $1$-skeletons of \Cech\
complexes, see][]{Penr:2003}.  Several results for the Delaunay
triangulation, some of which generalize to the Alpha complex, are reported in
\citep {Penr:Yuki:2001}. Regarding the support on the distribution of hypergraphs induced 
by the complexes, generally, this is an area of open research (personal communication with H. Edelsbrunner). 
Recent work by \cite{Kalh:2014} surveys some of this literature, focusing on random simplicial complexes. 
The monograph by \cite{Penr:2003} discusses the relationship between $\cQ$ and subgraph counts, 
the degree distribution, and the percolation threshold, in Chapters 3, 4 and 10, respectively
 
\citet [Chap.~3]{Penr:2003} gives conditions which guarantee the asymptotic
normality of the joint distribution of the numbers $Q_j$ of $j$-simplices
(edges, triads, \etc.), for \IID\ samples $\bV=(V_1,\ldots,V_\Num)$ from some
marginal distribution $\cQ_M$ on $\rr^\Dim$, as the number $\Num=|\bV|$ of vertices
grows and the radius $r_\Num$ shrinks.  

Simulation
studies suggest that the asymptotic results apply approximately for $\Num \ge
24$--$100$. By this we mean that sometimes $24$ is sufficient (the distribution of the 
vertices is approximately multivariate normal), and sometimes $100$ may be 
required (distribution of the vertices far from being multivariate normal).
%Theorems $3.4$, $3.5$, $3.12$ and $3.13$ of \citep %[\pg 52, 54, 65, 67]
%[Chap.~3] {Penr:2003} present features of the asymptotic distribution of
%subgraph counts for a RGG under the assumption that $\bV=(V_1,\ldots,V_\Num)$
%are \IID\ samples from some marginal distribution $\cQ_M$ on $\rr^\Dim$.
%Theorems $3.5$ and $3.13$ give conditions under which the asymptotic joint
%distribution of the subgraph counts for the numbers $Q_j$ of $j$-simplices is
%either multivariate normal or independent Poisson.  Simulation studies
%suggest that the asymptotic results of \citep {Penr:2003} apply approximately
%for vertex sets of cardinality at least $|\cV| \ge 24$--$100$.
%, depending on the values of various parameters.

\section{Geometric representations of random hypergraphs}\label{sec:StatModel} 
% Section 4
%\subsection{General setting}
We develop a Bayesian approach to the problem of inferring factorizations of 
distributions of the forms of \Eqn{eq:factors},
\begin{equation}\label{eq:theModel}\notag
f(x)= 
\prod_{a \in \sC(\cG)} \phi_a(x_a\mid \theta_a)
\text{\quad or\quad }
\frac{\prod_{a \in \sP(\cG)} \psi_a(x_a\mid \theta_a)}
        {\prod_{b \in \sS(\cG)} \psi_b(x_b\mid \theta_b)}.
\end{equation}
In each case we specify the prior density function as a product
\begin{equation}\label{eq:thePrior}
p(\theta,\cG) = p(\theta \mid \cG)\,p(\cG)
\end{equation}
of a conditional hyper Markov law for $\theta\in\Theta$ and a
marginal RGG law on $\cG$.  We use conventional methods to select the
specific hyper Markov distribution (hyper Inverse Wishart for multivariate
normal sampling distributions, for example) since our principal focus is on
prior distributions for the graphs, $p(\cG)$.  Every time we refers to hyper 
Markov  laws, it will be in the strong sense according to \citet {Dawi:Laur:1993}. We also present MCMC
algorithms for sampling from the posterior distribution on $\cG$, for
observed data.

\subsection{Prior specifications}
All the graphs in our statistical models are built from nerves constructed in
\Sec {s:rgg} from a random vertex set $\cV=\set{V_i}_{i=1}^\Num \subset
\rr^\Dim$ and radius $r>0$.  Since the nerve construction is invariant under
rigid transformations (this is, transformations that preserve angles as well as distances) of $\cV$ or simultaneous scale changes in $\cV$ and
$r$, restricting the support of the prior distribution on $\cV$ to the unit
ball $\BB \Dim$ does not reduce the model space:

\begin{prop}\label{p:feas}
  Every feasible graph in $\rr^\Dim$ may be represented in the
  form $\cG(\cV, r, \cA)$ for a collection $\cV$ of $\Num$ points in the unit
  ball $\BB \Dim$ and for $r=\ooNum$.
\end{prop}

\begin{proof}
Let $\cG=(\cV,\cE)\cong\cG(\cV, r, \cA)$ be a feasible graph with
$|\cV|=\Num$ vertices.  Every edge $\edg{v_i,v_j}\in\cE$ has length
$\dist{v_i,v_j}\le 2r$ so, by the triangle inequality, every connected
component $\Gamma_i$ of $\cG$ with $\Num_i$ vertices must have diameter no
greater than the longest possible path length, $2r(\Num_i-1)$, and so fits in
a ball $B_i$ of diameter $2r(\Num_i-1)$.  The union of these balls, centered
on a line segment and separated by $r(2+1/\Num)$, will have diameter less
than $r(2\Num-1)$.  By translation and linear rescaling we may take
$r=1/\Num$ and bound the diameter by $2$, completing the proof.
\end{proof}

We fix $r=\ooNum$ and simplify the notation by writing $\cG(\cV,\cA)$ instead
of $\cG(\cV,r,\cA)$ for $\cA=\Cech$ or $\cA=\Alpha$ or, if $\cA$ is
understood, simply $\cG(\cV)$.  Thus we can induce prior distributions on the
space of feasible graphs from distributions on configurations of $\Num$
points in the unit ball in $\rr^\Dim$.

For \IID\ uniform draws $\bV=(V_1,\dots,V_\Num)$ from $\BB\Dim$, the expected
number of edges of the graph $\cG(\bV,r,\cA)$ is bounded above by $\E
[\#\cE]\le\binom n2 (2r)^d$; for $r=\ooNum$ in dimension $\Dim=2$ this is
less than $\E [\#\cE]< 2$, leading to relatively sparse graphs.  We often
take larger values of $r$ (still small enough for empty graphs to be
feasible), to generate richer classes of graphs.  A limit to how large $r$
may be is given by the partial converse to \Prop{p:feas},

\begin{prop}\label{p:infeas}
  The empty graph on $\Num$ vertices cannot be expressed as $\cG(\cV, r,
  \Cech)$ for any $\cV\subset\BB \Dim$ with
  $r\ge\big(\Num^{1/\Dim}-1\big)^{-1}$.
 \end{prop}
\begin{proof} 
Let $\cV=\set{V_1,\dots,V_\Num}\subset\BB\Dim$ be a set of points and $r>0$ a
radius such that $\cG(\cV, r, \Cech)$ is the empty graph.  Then the balls
$\cV_i+r\BB\Dim$ are disjoint and their union
% $\cup[\cV_i+r\BB\Dim]$ 
with $\Dim$-dimensional volume $\Num \omega_\Dim r^\Dim$ lies wholly within
the ball $(1+r)\BB\Dim$ of volume $\omega_\Dim (1+r)^\Dim$ (where $\omega_d =
\pi^{\Dim/2} / \Gamma(1{+}\Dim/2)$ is the volume of the unit ball), so
$\Num <(1+\frac1r)^\Dim$.
\end{proof}
Slightly stronger, the empty graph may not be attained as $\cG(\cV, r,
\Cech)$ for any $r\ge1/[(\Num/p_\Dim)^{1/\Dim}-1]$ where $p_\Dim$ is the
maximum spherical packing density in $\rr^\Dim$.  For $\Dim=2$, this gives
the asymptotically sharp bound $r<1\Big/\Big[\sqrt{p\sqrt{12}/\pi}-1\Big]$.

\subsection{Sampling from prior and posterior distributions}\label{ss:samp}

Let $\cQ$ be a probability distribution on $\Num$-tuples in $\rr^\Dim$,
$p(\cG)$ the induced prior distribution on graphs
%for the \Cech\ complex construction 
$\cG(\bV,\Cech)$ for $\bV\sim\cQ$ with $r=\ooNum$, and let $p(\theta\mid\cG)$
be a conventional hyper Markov law (see below).  We wish to draw samples from
the prior distribution $p(\theta,\cG)$ of \Eqn {eq:thePrior} and from the
posterior distribution $p(\theta,\cG\mid\bx)$, given a vector $\bx=
(x_1,...,x_m)$ of \IID\ observations $x_j\iid f(x\mid\theta)$, using the
Metropolis\slash Hastings approach to MCMC \citep [Ch.~7] {Hast:1970,
  Robe:Case:2004}.

We begin with a random walk proposal distribution in $\BB \Dim$ starting at
an arbitrary point $v\in\BB \Dim$, that approximates the steps $\set{V\tim0,
V\tim1, V\tim2,...}$ of a diffusion $V\tim t$ on $\BB \Dim$ with uniform
stationary distribution and reflecting boundary conditions at the unit sphere
$\partial\BB \Dim$.

The random walk is conveniently parametrized in spherical coordinates with
radius $\rho\tim t=\parallel V\tim t \parallel$ and Euler angles--- in $\Dim{=}2$ dimensions,
angle $\varphi\tim t$--- at step $t$.  Informally, we take independent radial
random walk steps such that $(\rho\tim t)^\Dim$ is reflecting Brownian motion
on the unit interval (this ensures that the stationary distribution will be
$\Un(\BB \Dim)$) and, conditional on the radius, angular steps from Brownian
motion on the $\Dim$-sphere of radius $\rho\tim t$.

Fix some $\eta>0$.  In $\Dim=2$ dimensions the reflecting random walk proposal
$(\rho^*, \varphi^*)$ we used for step $(t+1)$, beginning at 
$(\rho\tim t,\varphi\tim t)$, is:
\[ 
  \rho^* = R\Big( [\rho\tim t]^2+\zeta_\rho\tim t\,\eta \Big)^{1/2},\qquad
  \varphi^* = \varphi\tim t + \zeta_\phi\tim t\, \eta/ \rho\tim t
  \] for \IID\ standard normal random variables $\set{\zeta_\rho\tim
  t,\zeta_\phi\tim t}$, where 
%\[R(x)=\left|x-2\left\lfloor\frac{x+1}2\right\rfloor\right|\]
\[R(x)=\left|x-2\left\lfloor\half(x+1)\right\rfloor\right|\]
%\[R(x)=\left|x-2\left\lfloor(x+1)/2\right\rfloor\right|\]
  is $x$ reflected (as many times as necessary) to the unit interval.
%sign is taken if $[\rho\tim t]^2+\eta \zeta_\rho\tim
%t\le1$, otherwise the negative sign.  
Similar expressions work in any dimension $\Dim$, with $\rho^* = R\Big(
  [\rho\tim t]^\Dim+\zeta_\rho\tim t\,\eta \Big)^{1/\Dim}$ and appropriate
step sizes for the $(\Dim-1)$ Euler angles.

\begin{figure}[!ht]
\begin{center}
  \includegraphics[height=110mm]{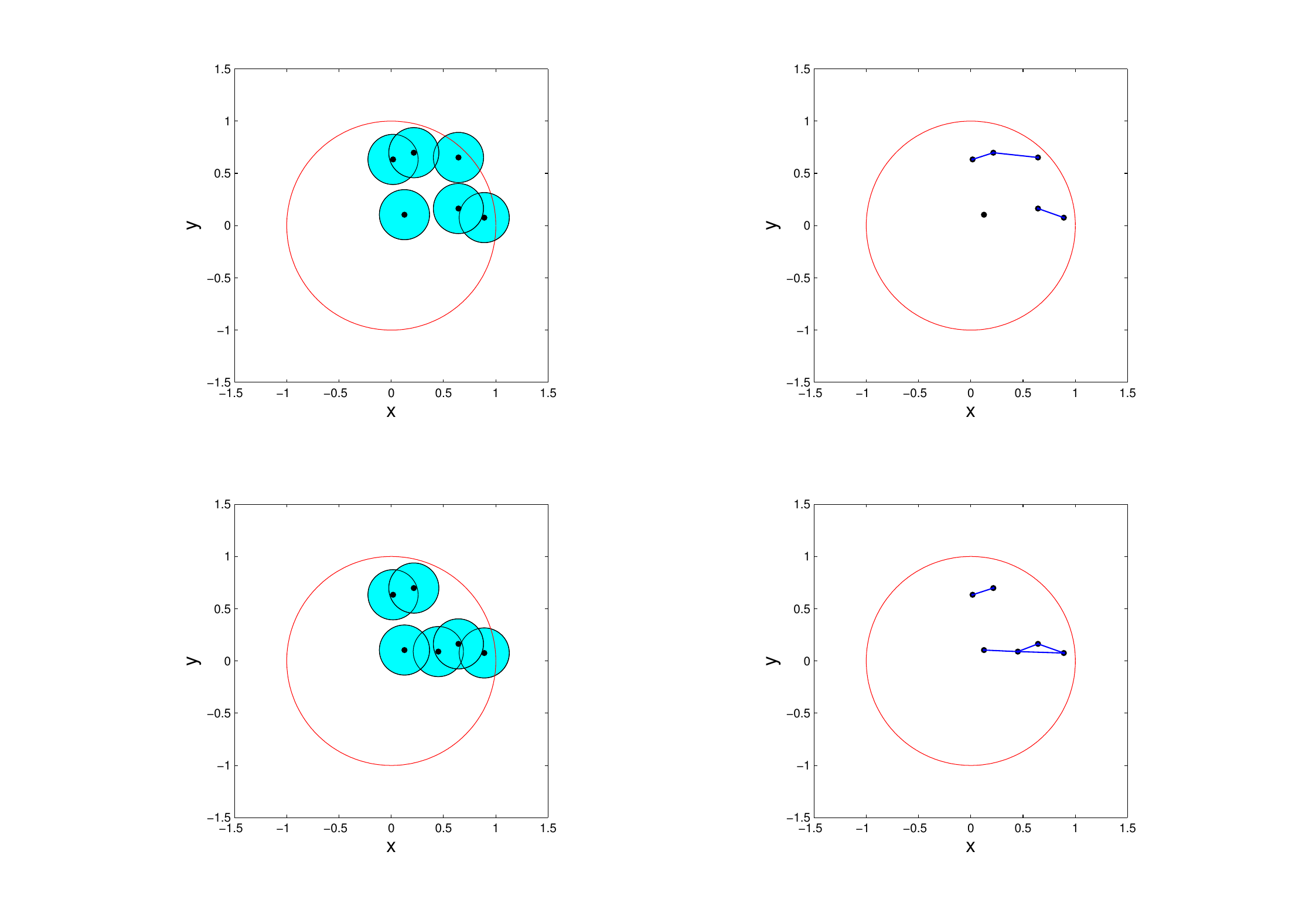}
  \caption{This figure illustrates a global move. (a) The current
    configuration of the points. (b) The graph implied by this 
    configuration. (c) The proposal configuration which is obtained by 
    randomly moving one vertex. (d) The graph implied by the proposed
    move.}\label{Fig:Global}
\end{center}
\end{figure}

For small $\eta>0$ this diffusion-inspired random walk generates local moves
under which the proposed new point $(\rho^*,\varphi^*)$ is quite close to
$(\rho\tim t,\varphi\tim t)$ with high probability. To help escape local
modes, and to simplify the proof of ergodicity below, we add the option of
more dramatic ``global'' moves by introducing at each time step a small
probability of replacing $(\rho\tim t,\varphi\tim t)$ with a random draw
$(\rho^*,\varphi^*)$ from the uniform distribution on $\BB \Dim$ (see
\Fig{Fig:Global}).  Let $q(\cV^*\mid \cV)$ denote the Lebesgue density at
$\cV^*\in(\BB \Dim)^\Num$ of one step of this hybrid random walk for
$\cV=(V_1,\dots,V_\Num)$, starting at $\cV\in(\BB \Dim)^\Num$.

\subsubsection{Prior sampling}\label{sss:priorSamp}

To draw sample graphs from the prior distribution begin with
$\bV\tim0\sim\cQ(d\bV)$ and, after each time step $t\ge0$, propose a new move
to $\bV^*\sim q(\bV^*\mid\bV\tim t)$.  The proposed move from $\bV\tim t$ (with
induced graph $\cG\tim t=\cG(\bV\tim t)$) to $\bV^*$ (and $\cG^*$) is
accepted (whereupon $\bV\tim {t+1}=\bV^*$) with probability $1\wedge H\tim
t$, the minimum of one and the Metropolis\slash Hastings ratio
\begin{equation}\label{e:MH-prior}\notag
  H\tim t=\frac{p(\bV^*)~~   q(\bV\tim t\mid\bV^*) } 
                {p(\bV\tim t)~q(\bV^*\mid\bV\tim t) }.
\end{equation} 
Otherwise $\bV\tim {t+1}=\bV\tim t$; in either case set $t\leftarrow t{+}1$
and repeat.  Note the proposal distribution $q(\cdot\mid\cdot)$ leaves the
uniform distribution invariant, so $H\tim t\equiv1$ for $\cQ(d\bV)\propto
d\bV$ and in that case every proposal is accepted.
% Some proposals are rejected for the repulsive processes, of course.

\subsubsection{Posterior sampling}\label{sss:postSamp}

After observing a random sample $X=\bx=(x_1,\dots,x_m)$ from the distribution
$x_j\sim f(x\mid \theta,\cG)$, let
\begin{align}
f(\bx\mid \theta,\cG)&=\prod_{i=1}^m f(x_i\mid \theta,\cG)\notag
\intertext{denote the likelihood function and}
  \cM(\cG)&=\int_{\Theta_\cG} f(\bx\mid \theta,\cG)\,
                      p(\theta \mid \cG)d\theta\label{e:MargLH}
\intertext{the marginal likelihood for $\cG$.  For posterior sampling of 
  graphs, a proposed move from $\bV\tim t$ to $\bV^*$ is accepted with
  probability $1\wedge H\tim t$ for} H\tim t&=\frac{\cM(\cG^*)~
  p(\bV^*)~~   q(\bV\tim t\mid\bV^*) } {\cM(\cG\tim t)~ p(\bV\tim
  t)~q(\bV^*\mid\bV\tim t) }.\label {e:MH-post}
\end{align} 
For multivariate normal data $X$ and hyper inverse Wishart hyper Markov law
$p(\theta\mid\cG)$, $\cM(\cG)$ from \Eqn{e:MargLH} can be expressed in closed
form for decomposable graphs $\cG(\bV)$.  efficient algorithms for evaluating
\Eqn{e:MargLH} are still available even if this condition fails.

The model %\Eqn{eq:theModel} 
will typically be of variable dimension, since
the parameter space $\Theta_\cG$ for the factors may depend on the graph
$\cG=\cG(\bV)$.  Not all proposed moves of the point configuration $\bV\tim
t\rightsquigarrow\bV^*$ will lead to a change in $\cG(\bV)$; for those that
do
 %in this case an ordinary Metropolis ratio of \Eqn{e:MH-post} is
 %sufficient.  For the case where the proposed move would change the graph and
 %hence the parameter space $\Theta_\cG$, a somewhat more involved
 %reversible-jump approach \citep {Gree:1995} is needed.  We 
we implement reversible-jump MCMC \citep {Gree:1995,Siss:2005} using the
auxiliary variable approach of \citet {Broo:Giud:Robe:2003} to simplify the
book-keeping needed for non-nested moves $\Theta_\cG\rightsquigarrow
\Theta_{\cG^*}$.

% end \hide

\subsection{Convergence of the Markov chain}

Denote by $\dot{\cG} (\Num,\Dim,\cA)$ the finite set of feasible graphs with
$\Num$ vertices in $\rr^\Dim$, \ie, those generated from $1$-skeletons of
$\cA$-complexes.
% isomorphic to the nerve of $\Num$ \emph{closed} balls in $\rr^\Dim$.
For each $\cG\in\dot{\cG} (\Num,\Dim,\cA)$ let $V_\cG\subset(\BB \Dim)^\Num$
denote the set of all points $\bV=\set{V_1,\dots,V_\Num}\in(\BB \Dim)^\Num$ for
which $\cG\cong \cG(\bV, \ooNum, \cA)$, and set $\mu(\cG)=\cQ\big(V_\cG\big)$.
%$\mu(\cG)=|V_\cG|\,\Gamma(1{+}\Dim/2)^\Num/\pi^{\Dim\Num/2}$
%
%where $|A|$ denote the $(\Dim{\times}\Num)$-dimensional Lebesgue measure of a
%set $A\subseteq(\BB \Dim)^\Num$.
%$\omega_\Num = \pi^{\Num/2} / \Gamma(1{+}\Num/2)$, so
%|(\BB \Dim)^\Num|= \pi^{\Dim\Num/2}/(\Gamma(1+m/2)^\Num)
Then
\begin{prop}\label{p:ergod}
  The sequence $\cG\tim t=\cG(\bV\tim t, \ooNum, \cA)$ induced by the prior
  MCMC procedure described in \Sec{sss:priorSamp} samples each feasible graph
  $\cG\in \dot{\cG} (\Num,\Dim,\cA)$ with asymptotic frequency $\mu(\cG)$.
  The posterior procedure described in \Sec{sss:postSamp} samples each
  feasible graph with asymptotic frequency $\mu(\cG\mid \bx)$, the posterior
  distribution of $\cG$ given the data $\bx$ and hyper Markov prior $p(\theta
  \mid\cG)$.
\end{prop}

\begin{proof}
Both statements follow from the Harris recurrence of the Markov chain
$\bV^{(t)}$ constructed in \Sec {ss:samp}.  For this it is enough to find a
strictly positive lower bound for the probability of transitioning from an
arbitrary point $V\in(\BB \Dim)^\Num$ to any open neighborhood of another
arbitrary point $V^*\in(\BB \Dim)^\Num$ \citep[Theorem 6.38, pg.  225]
{Robe:Case:2004}.  This follows immediately from our inclusion of the global
move in which all $\Num$ points $\set{V_i}$ are replaced with uniform draws from
$(\BB\Dim)^\Num$.
%
%Partitions of the state space are defined by the equivalence class of graphs
%they induce, $\bV_\cG = \set{\bV_j, \bV_k : \cG(\bV_j,r,\cA) \cong
%\cG(\bV_j,r,\cA)}$.  In this case the partitions are indexed over the set
%$\left\{\bV_\cG, \cG \in \dot{\cG} (\Num,\Dim,\cA) \right\}$.  Given the
%$\sigma$-finite measure $\mu$ on the configuration space is well defined the
%measure of each partition $\mu(\bV_\cG)$ is well defined.  Since the map from
%$\bV_\cG \rightarrow \cG$ is deterministic this induces the same measure on
%graphs, $\mu(\cG) = \mu(\bV_\cG)$. The result of this is that the probability
%of the algorithm visiting a graph $\cG$ is the asymptotic frequency
%$\mu(\cG)$.
\end{proof}

It is interesting to note that while the sequence $\cG\tim t=\cG(\bV\tim t,
\ooNum, \cA)$ is a hidden Markov process, it is not itself Markovian on the
finite state space $\dot{\cG} (\Num,\Dim,\cA)$; nevertheless it is ergodic,
by \Prop{p:ergod}.

\section{Results}\label{Sec:Simul}
%\section{Four Simulation Examples}\label{Sec:Simul}
% Section 5

Here we illustrate the use of the proposed parametrization using simulations and real data.
These numerical examples provide us with an opportunity to test priors that encourage sparsity, and MCMC algorithms that allow for local as well as
global moves by design.  

\subsection{Illustration of modeling advantages}

%\subsubsection{Example}
\subsubsection{The nerve determines the junction tree factorization}

Here we use a junction tree factorization with each univariate marginal $X_i$
associated to a point $V_i\in \rr^\Dim$ (the standard graphical models
approach).  In this case, specifying the class of sets to compute the nerve
and the value for $r$ determines a factorization for the joint density of
$\set{X_1,\ldots,X_\Num }$.  We illustrate with $\Num=5$ points in Euclidean space of
dimension $\Dim=2$.

Let $(X_1,X_2,X_3,X_4,X_5)\in\rr^2$ be a random vector with density $f(x)$
and consider the vertex set displayed in \Tab{tab:VertexExample} (also shown
as solid dots in \Figs {fig:Alpha_example} {fig:CechExample}).

\begin{figure}[htb]
\begin{center}
\includegraphics[height=80mm]{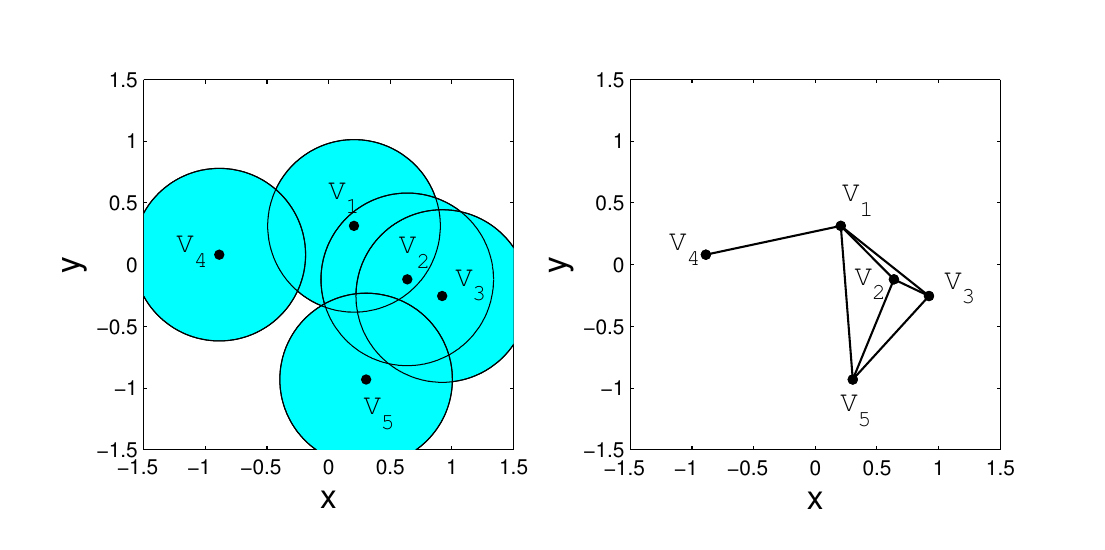}
\caption{(a) \Cech\ complex computed using $r=0.7$ and the
  vertex set listed in \Tab{tab:VertexExample}. (b) The $1$-Skeleton of the \Cech\ complex.}\label{fig:CechExample}
\end{center}
\end{figure}

For an Alpha complex with $r=0.5$ the junction tree factorization
\Eqn{eq:fac-jt} corresponding to the graph in \Fig {fig:Alpha_example} is
%\uhoh{Should these $f$'s be $\psi_a$'s as in (1.3)?}
\begin{displaymath}
%  f(x)=\frac{f(x_2,x_3,x_5)f(x_1,x_2) f(x_4)}{f(x_2)}.
  f(x)=\frac{\psi_{12}(x_1,x_2)
             \psi_{235}(x_2,x_3,x_5)
             \psi_{4}(x_4)
           }{\psi_{2}(x_2)},
\end{displaymath}
we will denote the factorization as $[1,2][2,3,5][4]$. In the case
where the factors are potential functions rather than marginals we
will use $\cs{\cdot}$ instead of $\cq{\cdot}$.
Similarly, for the \Cech\ complex and $r=0.7$ the factorization corresponding
to the graph in \Fig {fig:CechExample} is
\begin{displaymath}
%  f(x)=\frac{f(x_1,x_2,x_3,x_5)f(x_1,x_4)}{f(x_1)}.
  f(x)=\frac{\psi_{1235}(x_1,x_2,x_3,x_5)
             \psi_{14}(x_1,x_4)
           }{\psi_{1}(x_1)}.
\end{displaymath}

%\subsubsection{Simulation Study of Subgraph Counts for RGGs}
\subsubsection{Subgraph counts in RGGs are a function of $\cQ$}

% EDO : I AM ASSUMING THIS SECTION DOES NOT DESCRIBE NEW MATERIAL. I CHANGED THE TITLE AND SOME TEXT ACCORDINGLY

In this subsection we illustrate how the distribution of particular graph features changes as
a function of the sampling distribution of the random point set $\bV$ and
contrast this with \ER\ models.  Specifically we will focus on the number of
edges (2-cliques) $Q_2$ and the number of 3-cliques $Q_3$.
 
\begin{figure}[!ht]
\begin{center}
  \includegraphics[height=85mm]{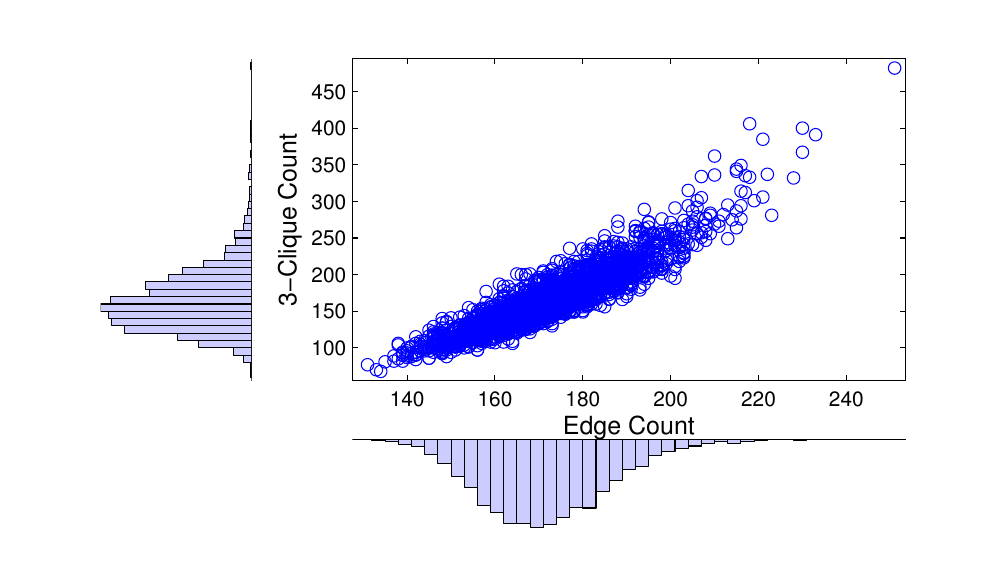}
  \caption{Edge counts and $3$-Clique counts from $2,500$ simulated samples
    of $\cG(\bV,1/\sqrt{2\cdot 75},\Cech)$ where
    $|\bV|=75$ and $V_i\iid\Un([0,1]^2)$, $1\leq i \leq 75$.
    The multivariate normal appears as a reasonable approximation for the
    joint distribution, as suggested by \citep[Thm. $3.13$] {Penr:2003}.
    \Cech\ radius is $r_n=1/\sqrt{2n}$.\label{fig:JointDist_Unif}}
\end{center}
\end{figure}

\begin{figure}[!ht]
\begin{center}
  \includegraphics[height=85mm]{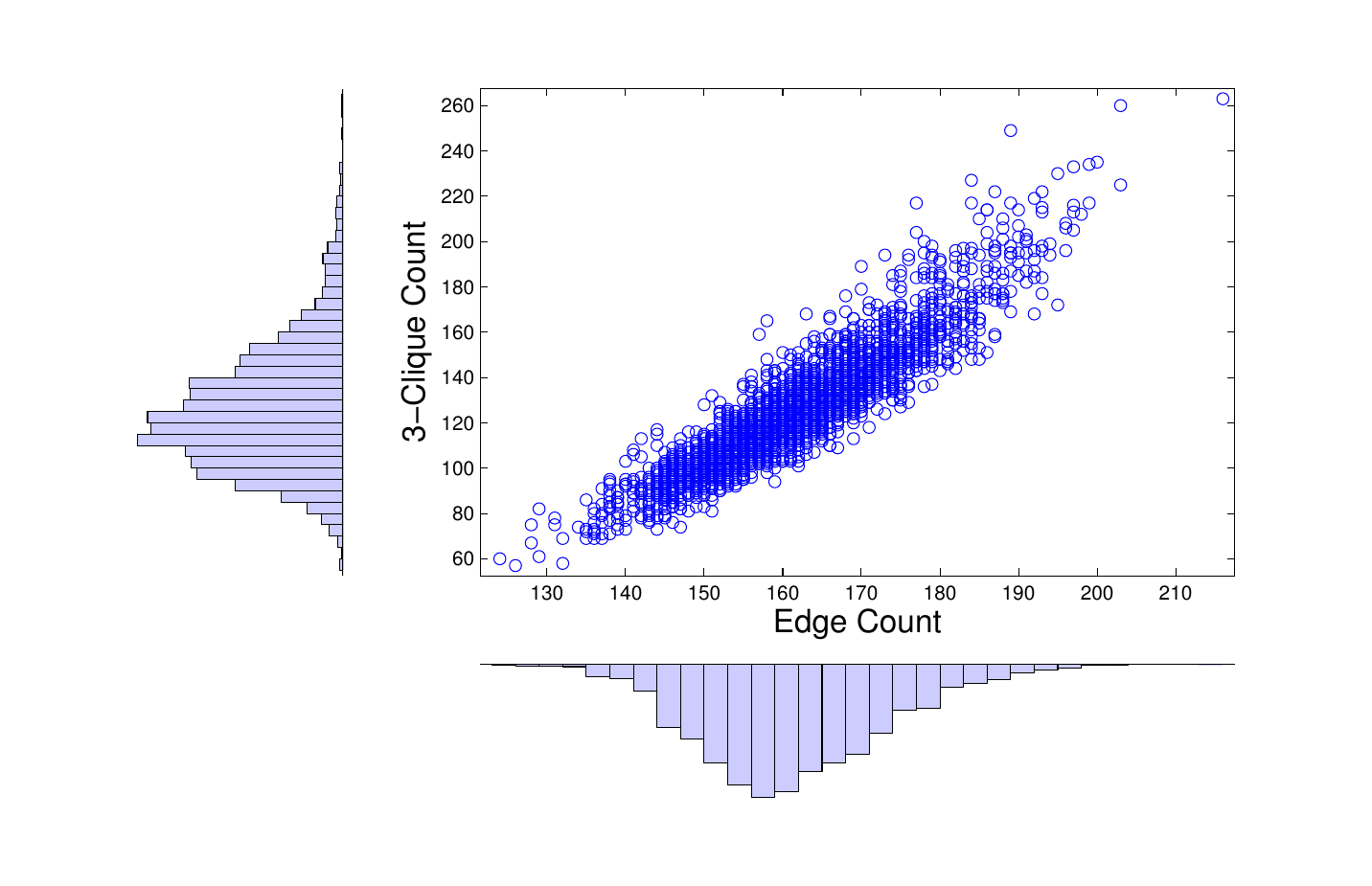}
  \caption{Edge counts and $3$-Clique counts from $2,500$ simulated samples
    of $\cG(\bV,1/\sqrt{2\cdot 75},\Cech)$ where
    $|\bV|=75$ and $V$ sampled from a Matt\'ern~III with hard-core radius
    $r=0.35$.\label{fig:JointDist_Matern}}
\end{center}
\end{figure}

\begin{figure}[t!]
\begin{center}
  \includegraphics[height=70mm]{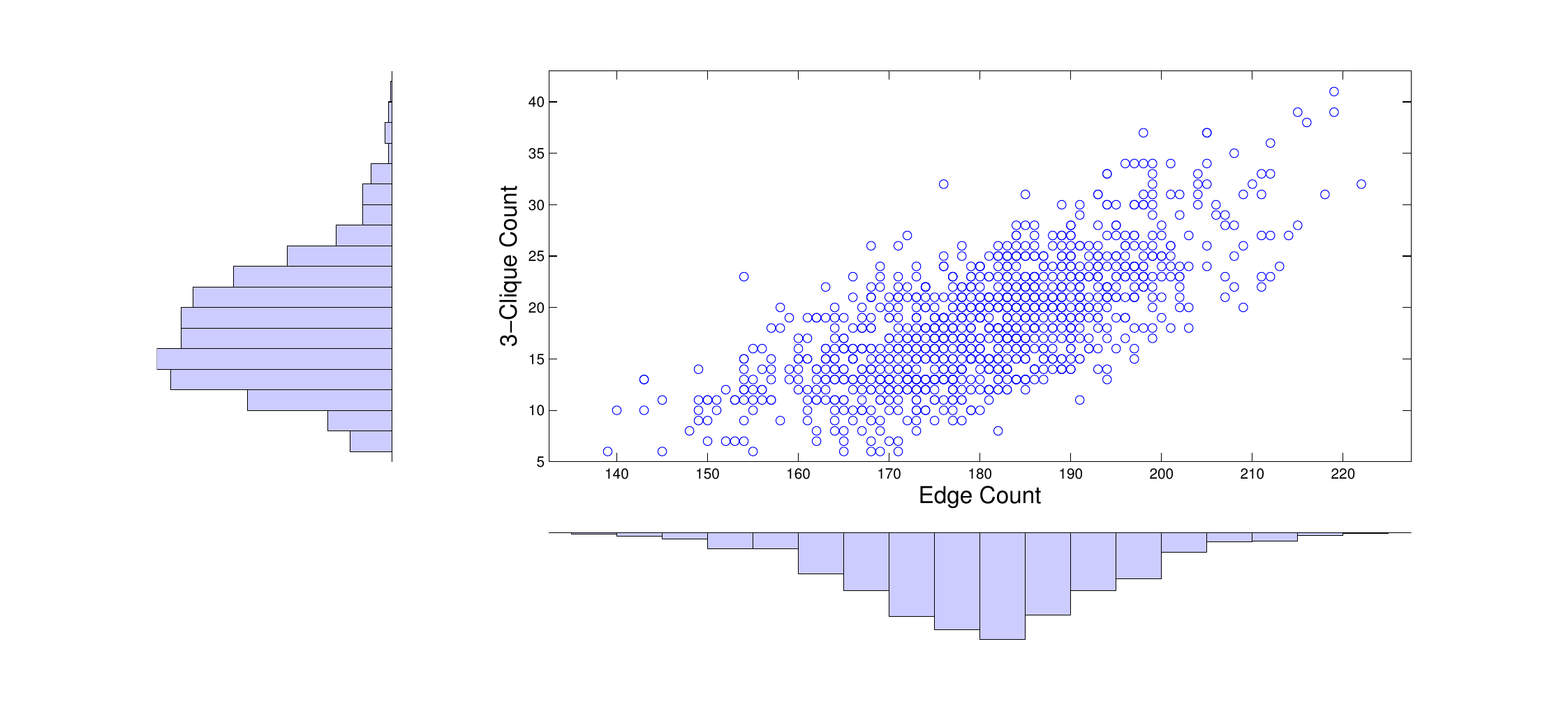}
  \caption{Edge counts and $3$-Clique counts from $1\,000$ simulated samples
    of an \ER\ graph with edge inclusion probability of
    $p=0.065$.\label{fig:JointDist_RE}}
\end{center}
\end{figure}

The two spatial processes we study for $\cQ$ are \IID\ uniform draws from the
unit square $[0,1]^2$ in the plane, and dependent draws from the Mat\'ern
type III hard-core repulsive process \citep {Hube:Wolp:2009}, using \Cech\
complexes with radius $r=1/\sqrt{150}\approx0.082$ in both cases to ensure
asymptotic normality \citep[Thm.\ 3.13] {Penr:2003}.  In our simulations we
vary both the number of variables (graph size) $\Num$ and the Mat\'ern III
hard core radius $\rho$.  Comparisons are made with an \ER\ model with a
common edge inclusion parameter.  \Tab{tabcompgraphstats} displays the
quartiles
%$25\th$, $50\th$,and $75\th$ percentiles 
for $Q_2$ and $Q_3$ as a function of the graph size $\Num$, hard core radius
$\rho$, and \ER\ edge inclusion probability $p$.  Figures
\ref{fig:JointDist_Unif}, \ref {fig:JointDist_Matern}, and
\ref{fig:JointDist_RE} show the joint distribution of $(Q_2,~Q_3)$ for
$\set{V_i}\iid \Un([0,1]^2)$, for a Mat\'ern III process with hard core radius
$\rho=0.35$, and for draws from an \ER\ model with inclusion probability
$\alpha=0.065$, respectively.

\clearpage

\begin{table}[t!]
\centering
\begin{tabular}{|r@{\thinspace}L |C| R R R | R R R |}
\hline
\multicolumn{2}{|c|}{\Strut Graph}&|\cV| &
\multicolumn{3}{c|}{Edges}&\multicolumn{3}{c|}{$3$-Cliques}\\
%&&&0.25 & 0.5 & 0.75 &0.25 & 0.5 & 0.75\\
&&&25\% & 50\% & 75\% &25\% & 50\% & 75\%\\
\hline\Strut
Uniform &         &75  &161&171 &182 &134 &160&190\\
Mat\'ern&(0.035)  &75  &154&161 &170 &110 &124&144\\
ER      &(0.050)  &75  &130&138 &146 &6   &8  &11\\
ER      &(0.065)  &75  &172&181 &189 &14  &18 &22\\
\hline\Strut
Uniform &         &50  &69 &75  &81  &34  &43 &57\\
Mat\'ern&(0.035)  &50  &66 &71  &76  &27  &35 &43\\
Mat\'ern&(0.050)  &50  &62 &67  &71  &22  &27 &33\\
ER      &(0.050)  &50  &56 &61  &67  &1   &2  &4\\
ER      &(0.065)  &50  &74 &79  &85  &3   &5  &7\\
\hline\Strut
Uniform &         &20  &9  &12  &14  &1   &2  &4\\
Mat\'ern&(0.035)  &20  &9  &11  &13  &1   &1  &3\\
Mat\'ern&(0.050)  &20  &8  &10  &12  &0   &1  &2\\
ER      &(0.050)  &20  &8  &9   &11  &0   &0  &0\\
ER      &(0.065)  &20  &10 &12  &15  &0   &0  &1\\
\hline
\multicolumn{8}{c}{}
\end{tabular}
%\end{center}\par
\caption{Summaries of the empirical distribution of edge and $3$-clique
  counts for \Cech\ complex random geometric graphs with radius $r=0.082$,
  for vertex sets sampled from \IID\ draws from the unit square from: a
  uniform distribution, a hard core process with radius $\rho=0.035$, and
  from \ER\ (ER) with common edge inclusion probabilities of $\alpha=0.050$ and
  $\alpha=0.065$.\label{tabcompgraphstats}}
\end{table}

These simulations illustrate that by varying the distribution $\cQ$ we can
control the joint distribution of graph features.  The repulsive and \IID\
uniform distribution have very similar edge distributions, for example (see
\Figs {fig:JointDist_Unif} {fig:JointDist_Matern}), while (as anticipated)
the repulsive process penalizes large cliques.  Joint control of these
features is not possible with an \ER\ model with a common edge inclusion
probability and it is not obvious how to encode this type of information in
the concordance function approach of \citet {Mukh:Spee:2008}.

In \Sec{Sec:Filtration} we proposed a procedure for generating decomposable
graphs, and noted that the graphs induced by this algorithm are similar to
those constructed without the decomposability restriction.  In \Fig
{fig:Dist_Edges} we display a simulation study of the distribution of edge
counts for a RGG and the restriction to decomposable graphs.  These
distributions are very similar.

%\subsection{Numerical illustrations}
\subsection{Simulation studies}

We develop four examples.  The first example illustrates that
our method works when the graph encoding the Markov structure of underlying
density is contained in the space of graphs spanned by the nerve used to fit
the model.  In the second example we apply our method to Gaussian Graphical
Models.  The third example shows that the nerve hypergraph (not just the
$1$-skeleton) can be used to induce different groupings in the terms of a
factorization, and therefore a way to encode dependence features that go
beyond pairwise relationships.  In the fourth example we compare results
obtained by using different filtrations to induce priors over different
spaces of graphs.

\subsubsection{$\cG$ is in the Space Generated by $\cA$}\label
{Sec:ExampleIlust}

Let $\left(X_1,\ldots,X_{10} \right)$ be a random vector whose distribution
has factorization: 
\begin{subequations}\label{graph1}
\begin{equation}\label{eq:truemod1}
  f_\theta(\bx)=
%    \frac{\psi_\theta(x_1,x_3,x_{10})
%    \psi_\theta(x_1,x_8,x_{10})\psi_\theta(x_2, x_7)
%    \psi_\theta(x_2,x_4,x_6)\psi_\theta(x_5,x_9)}
%   {\psi_\theta(x_1,x_{10})\psi_\theta(x_2)}
\frac{\psi_\theta(x_1,x_4,x_{10})  \psi_\theta(x_1,x_8,x_{10})
\psi_\theta(x_4,x_5) \psi_\theta(x_8,x_9) \psi_\theta(x_2,x_3,x_9)
\psi_\theta(x_6) \psi_\theta(x_7)  }{\psi_\theta(x_4) \psi_\theta(x_8)
\psi_\theta(x_9) \psi_\theta(x_1,x_{10}) }
\end{equation}

\begin{figure}[htb]
\begin{center}
    \includegraphics[height=120mm]{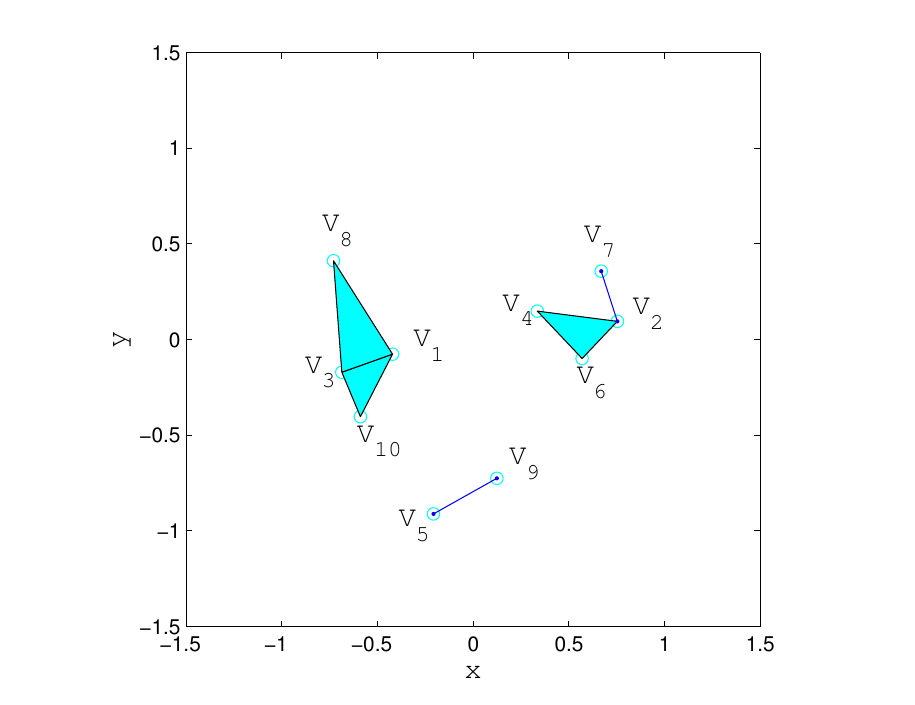}
    \caption{Geometric graph representing the model given in
      \Eqn{eq:truemod1}.  For this example graphs are constructed to be
      decomposable and the clique marginals are specified as Clayton
      copulas.}\label{fig:Toy2True}
\end{center}
\end{figure}

The Markov structure of \Eqn{eq:truemod1} can be encoded by the geometric
graph displayed in \Fig{fig:Toy2True}.  We transform variables if necessary
to achieve standard $\Un(0,1)$ marginal distributions for each $X_i$, and
model clique joint marginals with a Clayton copula \citetext {\citealp
  {Clay:1978}, or \citealp{Nels:1999}, \S4.6}, the exchangeable multivariate
model with joint distribution function
\begin{align}
 \Psi_\theta(\bx_I) &=
   \left(1-\Num_I+\sum_{i\in I} x_i^{-\theta}\right)^{-1/\theta}\notag
\intertext{and density function}
 \psi_\theta(\bx_I) &=
\iffalse
 \frac{\theta^{\Num_I}~\Gamma(\Num_I+1/\theta)/\Gamma(1/\theta)}
      {\big(1-\Num_I+\sum_{i\in I} x_i^{-\theta}\big)^{\Num_I+1/\theta}}
 ~ \prod_{i\in I} x_i{} ^{-1-\theta}
\else
 \theta^{\Num_I}~
 \frac{\Gamma(\Num_I+1/\theta)}
      {\Gamma(1/\theta)}
   \left(1-\Num_I+\sum_{i\in I} x_i^{-\theta}\right)^{-\Num_I-1/\theta}
 ~ \left(\prod_{i\in I} x_i\right) ^{-1-\theta}\label{e:Clay}
\fi
\end{align}
\end{subequations}
on $[0,1]^{\Num_I}$ for some $\theta\in\Theta=(0,\infty)$, for each clique
$\cq{v_i:~i\in I}$ of size $\Num_I$.

\begin{figure}[ht]
\begin{center} % was: {Talk_10v_samples.eps}
 \includegraphics[height=150mm]{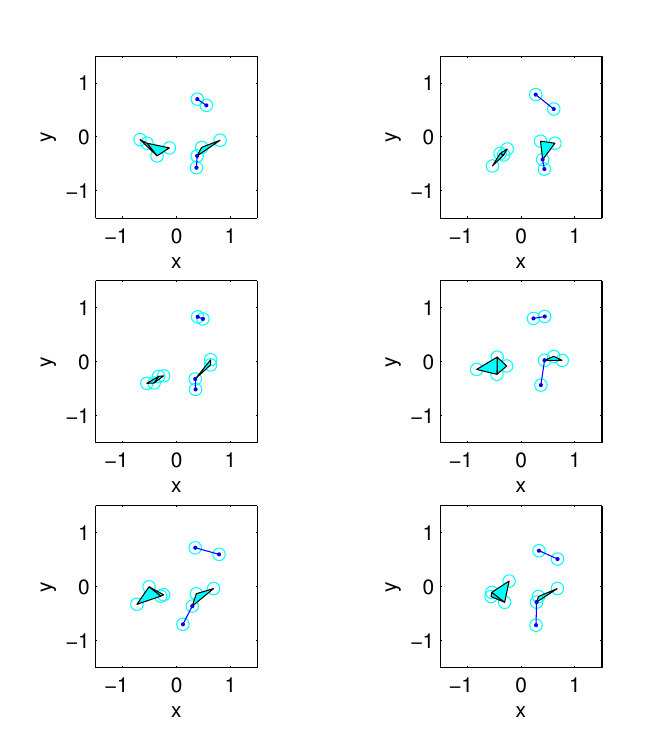}
  \caption{Geometric graphs corresponding to snapshots of posterior samples (one every 100 iterations) from model of \Eqn{eq:truemod1}. For this example graphs are constructed to be decomposable and the clique marginals are specified as Clayton copulas.}\label{fig:Toy3Samples1}
\end{center}
\end{figure}

We drew $250$ samples from model \Eqn{graph1} with $\theta =4$.  For
inference about $\cG$ we set $\cA=\Alpha$ and $r=0.30$, with independent
uniform prior distributions for the vertices $V_i\iid \Un(\BB2)$ on the unit
ball in the plane.  We used the random walk described in \Sec {ss:samp} to
draw posterior samples with \Alg {filtrationalg} applied to enforce
decomposability.  To estimate $\theta$ we take a unit Exponential prior
distribution $\theta\sim\Ex(1)$ and employ a Metropolis\slash Hastings
approach using a symmetric random walk proposal distribution with reflecting
boundary conditions at $\theta=0$,
\[ \theta^* = \big| \theta\tim t + \varepsilon\big|,\] with
$\varepsilon_{t}\sim \Un(-\beta,\beta)$ for fixed $\beta >0$.  We drew
$1\,000$ samples after a burn-in period of $25\,000$ draws.  The three models
with the highest posterior probabilities are displayed in
\Tab{tab:TopTopol3}.  The geometric graphs computed from six posterior
samples (one every $100$ draws) are shown in \Fig {fig:Toy3Samples1}; note
that the computed nerves appear to stabilize after a few hundred iterations
while the actual position of the vertex set continues to vary.

\begin{table}[ht]
\begin{center}
\begin{tabular}{|L C|}
\hline
\Strut\text{Graph Topology}&\text{Posterior Probability}\\
\hline\Strut
% \CQ{1,3,10}\CQ{1,3,8}\CQ{2,4,6}\CQ{2,7}\CQ{5,9}& 0.964\\
% \cq{1,3,10}\cq{1,3,8}\cq{2,4,6}\cq{2,4,7}\cq{5,9}  & 0.017\\
% \cq{1,3,10}\cq{1,3,8}\cq{2,4,6}\cq{2,7}\cq{5}\cq{9}& 0.015\\
%%%% [1,10,4][1,10,8][4,5][8,9][2,3,9][6][7]     0.963
%%%% [1,10,4][1,10,8][4,5][8,9][2,3,9][6][5,7]   0.021
%%%% [1,10,4][1,8][4,5][8,9][2,3,9][6][7]        0.010
 \CQ{1,4,10}\CQ{1,8,10}\CQ{4,5}\CQ{8,9}\CQ{2,3,9}\CQ{6}\CQ{7}&0.963\\
 \cq{1,4,10}\cq{1,8,10}\cq{4,5}\cq{8,9}\cq{2,3,9}\cq{6}\cq{5,7}&0.021\\
 \cq{1,4,10}\cq{1,8}\cq{4,5}\cq{8,9}\cq{2,3,9}\cq{6}\cq{7}& 0.010\\
\hline
\multicolumn{2}{c}{}
\end{tabular}
\caption{The three models with highest estimated posterior probability.  The
  true model is shown in bold (see \Fig{fig:Toy2True}).  Here
  $\theta=4$.}\label{tab:TopTopol3}\notag
\end{center}
\end{table}

\subsubsection{Gaussian graphical model}\label{Sec:Gaussian}

We use our procedure to perform model selection for the Gaussian graphical
model $X\sim\No(0,\Sigma_{\cG})$, where $\cG$ encodes the zeros in
$\Sigma^{-1}$.  We adopt a Hyper Inverse Wishart (HIW) prior distribution for
$\Sigma\mid\cG$.  The marginal likelihood
% (in the parametrization of \citealt[Eqn~(12)] {Atay:Mass:2005})
\citep [in the parametrization of] [Eqn~(12)] {Atay:Mass:2005}
is given by
\begin{align}
   \cM(\bV) &= (2\pi)^{-\Num N/2} ~ \frac{I_{\cG(\bV)}
   (\delta +N, D+ X^TX)}{I_{\cG(\bV)}(\delta ,D)},\label{eq:GaussMargLik}
\intertext{where}
  I_\cG(\delta ,D) &= \int_{M^+(\cG)} |\Sigma|^{(\delta-2)/2} e^{-\frac12<\Sigma,D>}
                   \,d\Sigma\notag
\end{align} 
denotes the HIW normalizing constant.  This quantity is available in closed
form for weakly decomposable graphs $\cG(\bV)$, but for our unrestricted
graphs \Eqn {eq:GaussMargLik} must be approximated via simulation.  For our
low-dimensional examples the method of \citep {Atay:Mass:2005} suffices; for
larger numbers of variables we recommend that of \citep{Carv:Mass:West:2007}.
We set $\delta=3$ and $D=0.4I_6+0.6J_6$ ($I_6$ and $J_6$ denote the identity
matrix and the matrix of all ones, respectively).

\begin{figure}[htb]
\begin{center}
  \includegraphics[height=75mm]{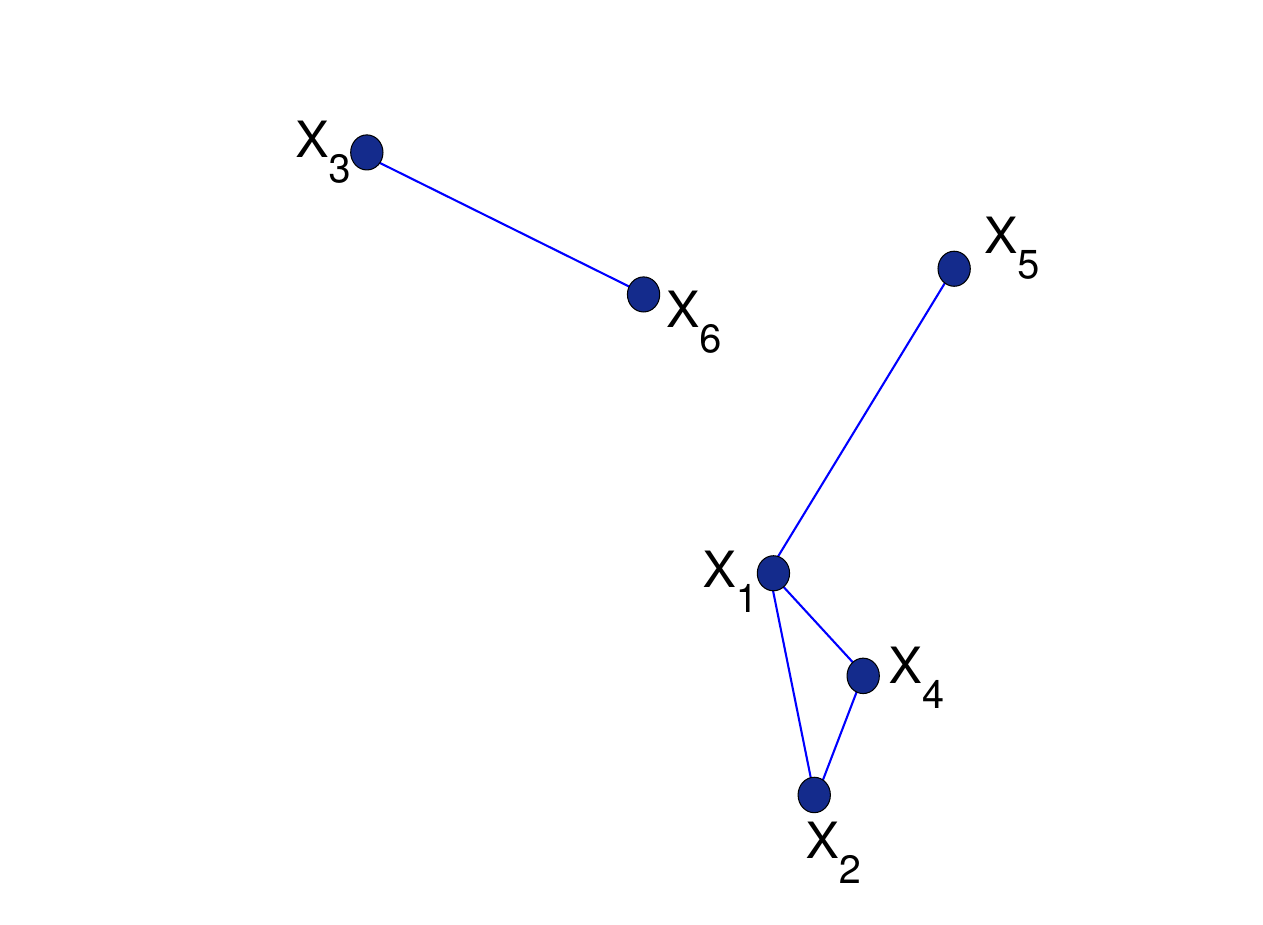}
  \caption{Graph encoding the Markov structure of the model given by
    precision matrix \Eqn{e:prec}.}\label{fig:GaussGraph}
\end{center}
\end{figure}

We sampled $300$ observations from a Multivariate Normal with mean zero and
precision matrix
\begin{equation}\label{e:prec}
\begin{pmatrix}
18.18 &  -6.55 & 0     & 2.26 & -6.27 &  0 \\
-6.55 &  14.21 & 0     & -4.90& 0     &  0 \\ 
0     &    0   & 10.47 & 0    & 0     & -3.65\\
2.26  & -4.90  & 0     & 10.69& 0     &  0 \\
-6.27 & 0      & 0     & 0    &27.26  &  0 \\
0     & 0      &-3.65  & 0    & 0     &7.41
\end{pmatrix}
\end{equation}
whose conditional independence structure is given by the graph in \Fig
{fig:GaussGraph}.  We fit the model described in \Sec{sec:StatModel} using a
uniform prior for each $V_i\in\BB2$ and $r=0.25$.  We employed hybrid random
walk proposals in which we move all five vertices $\set{V_i}$ independently
according to the diffusion-inspired random walk described in \Sec{ss:samp}
with probability $0.85$; replace one uniformly selected vertex $V_i$ with a
uniform draw from $\Un(\BB2)$ with probability $0.05$; and replace all five
vertices with independent unoform draws from $\Un(\BB2)$ with probability
$0.10$.  We sampled $1\,000$ observations from the posterior after a burn in
of $750\,000$.  Results are summarized in \Tab{tab:Gauss1}

\begin{table}[ht]
\begin{center}
\begin{tabular}{|L C|}
\hline
\text{Graph Topology}&\text{Posterior Probability}\\
\hline\Strut
 \CQ{1,2,4}\CQ{1,5}\CQ{3,6} & 0.152 \\
 \cq{1,5}\cq{2,3,4}\cq{2,3,6} & 0.072 \\
 \cq{1,2,3,4,6}\cq{1,5} & 0.069 \\
 \cq{1,4}\cq{2,4}\cq{2,3,6} & 0.055 \\
 \cq{1,2,4}\cq{2,3,4}\cq{1,5}\cq{3,6} & 0.052 \\
\hline\multicolumn{2}{c}{}
\end{tabular}
\caption{The five models with highest estimated posterior probability.
% In this case the true model is $\cq{1,2,4}\cq{1,5}\cq{3,6}$.
The true model is shown in bold.\label{tab:Gauss1}}
\end{center}
\end{table}
 
\subsubsection{Factorization Based on Nerves}\label{Sec:FacNerve}
While Gaussian joint distributions are determined entirely by the bivariate
marginals, and so only edges appear in their complete-set factorizations (see
\Eqn{eq:fac-comp}); more complex dependencies are possible for other
distributions.  The familiar example of the joint distribution of three
Bernoulli variables $X_1$, $X_2$, $X_3$ each with mean $1/2$, with $X_1$ and
$X_2$ independent but $X_3=(X_1-X_2)^2$ (so that $\set{X_i}$ are only
\emph{pairwise} independent) has only the complete set $\cs{1,2,3}$ in its
factorization.

\begin{figure}[b!]
\begin{center}
  \includegraphics[height=70mm]{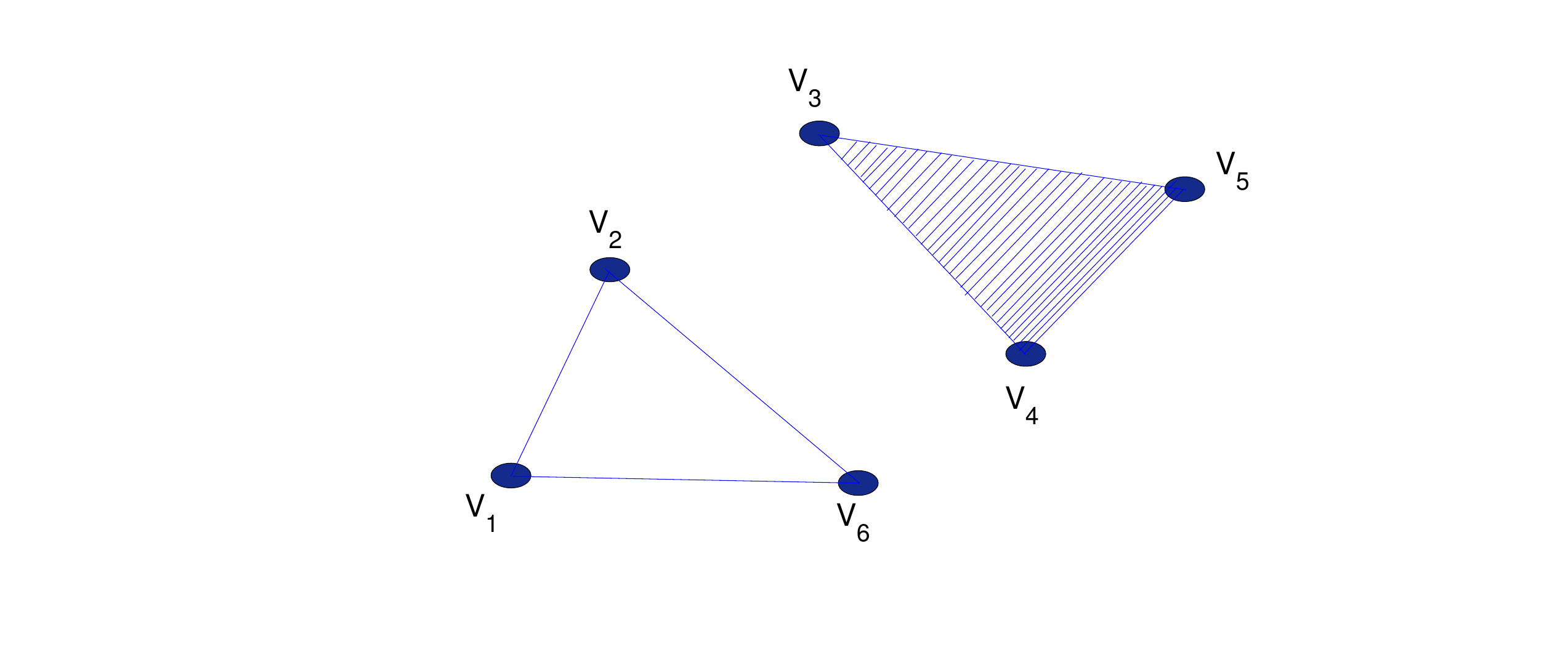}
  \caption{Graph encoding the Markov structure of the model given in
     \Eqn{e:HyperClay}.}\label{fig:GraphEx}
\end{center}
\end{figure}

Consider now a model with the graphical structure illustrated in \Fig
{fig:GraphEx} whose density function, if it is continuous and strictly
positive (\citep[see][Prop.~$3.1$]{Laur:1996}), admits the complete-set
factorization:
\begin{subequations}\label{e:HyperClay}
\begin{equation}\label{Eq:ModelHyper}
  f(x\mid\cG,\theta)= c_\cG\,\phi(x_1,x_2)
                             \phi(x_1,x_6)
                             \phi(x_2,x_6)\cdot
                             \phi(x_3,x_4,x_5).
\end{equation}
%The marginal distribution of $\left( X_3, X_4, X_5 \right)$ cannot be
%Gaussian, but that of $\left( X_1, X_2, X_6 \right)$ could.  
For illustration we will take each $\phi(\cdot)$ to be a Clayton copula
density (see \Eqn{e:Clay}).  %\citep {Clay:1978} and \pg152 in \citep {Nels:1999}).
For simplicity we specify the same value $\theta=4$ for each association
parameter, so $f(x\mid\cG,\theta)$ is given by \Eqn{Eq:ModelHyper} with
\iffalse % Note: Just here for reference, shouldn't appear in paper
\begin{align*}
C(u_1,...,u_n) &= \left(1-n+\sum u_j{}^{-\theta}\right)^{-1/\theta}\\
\phi(u_1,...,u_n) &= \frac{\partial^n}{\partial u_1\cdots\partial u_n}
                      C(\vec u)\\
                  &=\frac{\Gamma(\frac1\theta+n)\,\theta^n}
                         {\Gamma(\frac1\theta)}
                    \left(\prod u_j\right)^{-\theta-1}
                    \left(1-n+\sum u_j{}^{-\theta}\right)^{-n-1/\theta}
\end{align*}
\fi
\begin{alignat}5
\phi(x,y)   &=  5&&(x^{-4}+y^{-4}-1)^{-9/4}&&(x\,y)^{-5} \label{e:Clay2}\\
\phi(x,y,z) &= 30&&(x^{-4}+y^{-4}+z^{-4}-2)^{-13/4}&&(x\,y\,z)^{-5}.\label{e:Clay3}
\end{alignat}
\end{subequations}
%%f(x\mid\cG,\theta) & = c_\cG\, C^3(x_3,x_4,x_5)^{1+3\theta}\left( C^2(x_2,x_6)
%%  C^2(x_1,x_2)C^2(x_1,x_6) \right)^{1+2\theta} \nonumber \\
%%     & \times \left(x_1,x_2,x_5,x_1^2,x_2^2,x_6^2\right)
%%              ^{-(1+\theta)}(1+\theta)^4(1+2\theta).
%%\end{align}

\begin{figure}[t!]
\begin{center}
  \includegraphics[height=100mm]{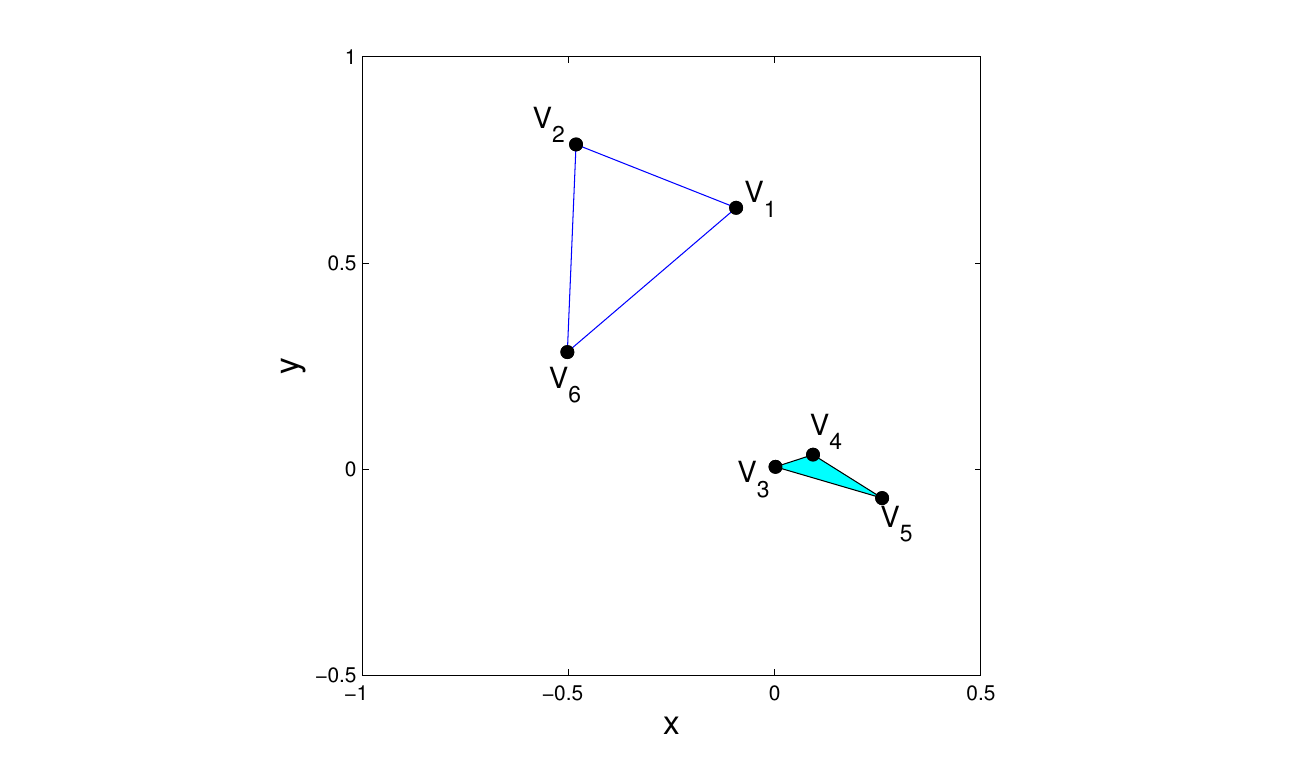}
  \caption{Alpha complex corresponding to the vertex set in \Tab
    {tab:VertexNervEx} and $r=\sqrt{0.075}$.}\label{fig:NervEx}
\end{center}
\end{figure}

%The challenge now is to infer a factorization over complete sets of a graph
%in the case where grouping the factors differently may lead to a different
%functional form for the density.  
In earlier examples we associated graphical structures (\ie, edge sets) with
$1$-skeletons of nerves.  We now associate \emph{hyper\/}graphical structures
(\ie, abstract simplicial complexes that may include higher-order simplexes)
with the entire nerves, with maximal simplices associated with complete-set
factors.  For example: the Alpha complex computed from the vertex set
displayed in \Tab {tab:VertexNervEx} with $r=0.40$ has $\cs{3,4,5} \cs{1,2}
\cs{1,6} \cs{2,6}$ as its maximal simplices (\Fig{fig:NervEx}).  By
associating a Clayton copula to each of these hyperedges we recover the model
shown in \Eqn{e:HyperClay}.

\begin{table}[t!]
\begin{center}
\begin{tabular}{|c|RRRRRR|}
\hline\Strut
Coordinate & V_1\quad{} & V_2\quad{} & V_3\quad{}
           & V_4\quad{} & V_5\quad{} & V_6\quad{} \\
\hline\Strut
 $x$ & -0.0936 &  -0.4817 &  0.0019 &  0.0930  &  0.2605 & -0.5028 \\
 $y$ &  0.6340 &   0.7876 &  0.0055 &  0.0351  & -0.0702 &  0.2839 \\
\hline\multicolumn{7}{c}{}
\end{tabular}
\caption{Vertex set used for generating a factorization based on
  nerves.}\label{tab:VertexNervEx}
\end{center}
\end{table}

We use the same prior and proposal distributions constructed in \Sec
{ss:samp} from point distributions in $\rr^\Dim$; what has changed is the way
the nerve is being used: as a hypergraph whose maximal hyperedges represent
factors.  One complicating factor is the need to evaluate the normalizing
factor $c_\cG$ for each graph $\cG$ we encounter during the simulation;
unavailable in closed form, we use Monte Carlo importance sampling to
evaluate $c_\cG$ for each new graph, and store the result to be reused when
$\cG$ recurs.

We anticipate that uniform draws $V_i\iid\Un(\BB2)$ will give high
probability to clusters of three or more points within a ball of radius $r$,
favoring higher-dimensional features (triangles and tetrahedra) in the
induced hypergraph encoding the Markov structure of $\set{X_i}$.  To explore
this phenomenon, we compare results for uniform draws with those from a
repulsive process under which clusters of three or more points are unlikely
to lie within a ball of radius $r$, hence favoring hypergraphs with only
edges.

We began by sampling $650$ observations from model \Eqn{e:HyperClay} with
$\cA=\Alpha$ and $r=0.40$, with independent uniform prior distributions for
the vertices $V_i\iid \Un(\BB2)$ on the unit ball in the plane.  The
Metropolis\slash Hastings proposals for the vertex set are given by a mixture
scheme:
\begin{itemize}
\item A random walk for each $V_i$ as described in \Sec {ss:samp}, with step
  size $\eta=0.020$.  This proposal is picked with probability $0.94$
\item An integer $1\le k\le 6$ is chosen uniformly and, given $k$, a subset
  of size $k$ from $\set{1,2,3,4,5,6}$ is sampled uniformly; the vertices
  corresponding to those indices are replaced with random independent draws
  from $\Un(\BB2)$.  This proposal is picked with probability $0.06$, $0.01$
  for each $k$.
\end{itemize}
For $\theta$ we used the same standard exponential prior distribution and
reflecting uniform random walk proposals described in Example \ref
{Sec:ExampleIlust}.

Using $5\,000$ posterior samples after a burn-in period of $95\,000$
iterations, the models with highest posterior probability are summarized in
\Tab{tab:FactUnif}.

\begin{table}[b!]
\begin{center}
\begin{tabular}{|L C|}
\hline\Strut
\text{Maximal Simplices}&\text{Posterior Probability}\\
\hline\Strut
%%%%{3,4,5}{1,2}{2,6}{1,6}       0.6088
%%%%{1,2,6}{3,4}{4,5}{3,5}       0.1612
%%%%{3,5}{1,6}{3,4}{1,2}{2,6}    0.1370
 \CS{3,4,5}\CS{1,2}\CS{2,6}\CS{1,6}       & 0.609\\
 \cs{1,2,6}\cs{3,4}\cs{4,5}\cs{3,5}       & 0.161\\
 \cs{3,5}\cs{1,6}\cs{3,4}\cs{1,2}\cs{2,6} & 0.137\\
%%%% \CS{3,4,5}\CS{1,2}\CS{2,6}\CS{1,6} & 0.374 \\
%%%% \cs{3,4,5}\cs{1,2,6} & 0.259 \\
%%%% \cs{3,4,5}\cs{1,6}\cs{1,2} & 0.068 \\
%%%% \cs{1,2,6}\cs{4,5}\cs{3,5}\cs{3,4} & 0.042 \\
\hline
\multicolumn{2}{c}{}
\end{tabular}
\caption{Highest posterior factorizations with uniform prior
 for model of \Eqn{e:HyperClay} and \Fig{fig:GraphEx}
 (true model is shown in bold).}\label{tab:FactUnif}
\end{center}
\end{table}

To penalize higher-order simplexes we used a Strauss repulsive process \citep
{Stra:1975} conditioned to have $\Num$ points in $\BB \Dim$ as prior
distribution for the vertex set, with Lebesgue density
\begin{displaymath}
  g(v)\propto \gamma^{\#\set{(i,j):~ \dist{v_i,v_j}<2R }}
\end{displaymath}
for some $0<\gamma\le1$, penalizing each pair of points closer than $2R$.
Simulation results for this prior (with $R=0.7r$ and $\gamma=0.75$) are
summarized in \Tab {tab:FactStrauss}.  The posterior mode is far more distinct
for this prior than for the uniform prior shown in \Tab {tab:FactUnif}.

\begin{table}[t!]
\begin{center}
\begin{tabular}{|L C|}
\hline\Strut
\text{Maximal Simplices}&\text{Posterior Probability}\\
\hline\Strut
%%{3,4,5}{1,2}{2,6}{1,6}     0.8237
%%{1,2,6}{3,4,5}             0.1114
%%{1,2,6}{3,4}{4,5}{3,5}     0.0024
 \CS{3,4,5}\CS{1,2}\CS{2,6}\CS{1,6} & 0.824 \\
 \cs{1,2,6}\cs{3,4,5}               & 0.111 \\
 \cs{1,2,6}\cs{3,4}\cs{3,5}\cs{4,5} & 0.002 \\
%%%% \CS{3,4,5}\CS{1,2}\CS{2,6}\CS{1,6} & 0.368 \\
%%%% \cs{3,4,5}\cs{1,2,6} & 0.192 \\
%%%% \cs{1,2,6}\cs{3,4}\cs{3,5}\cs{4,5} & 0.183 \\
\hline\multicolumn{2}{c}{}
\end{tabular}
\caption{Highest posterior factorizations with Strauss prior
 (true model is shown in bold).}\label{tab:FactStrauss}
\end{center}
\end{table}
In a further experiment with $\gamma=0.35$, the posterior was concentrated on 
factorizations without any triads.
% (\ie, a heavier penalty on triads), 
%in the factorization given by the complete sets $\cs{1,2}$, $\cs{1,6}$,
%$\cs{2,6}$, $\cs{3,4}$, $\cs{3,5}$, and $\cs{4,5}$.

\subsubsection{$\cG$ Outside the Space Generated by $\cA$}\label
{Sec:3Filtrations}
In the simulation studies of Sections
\ref{Sec:ExampleIlust}--\ref{Sec:FacNerve} the class of sets $\cA$ used to
compute the nerve was known.  In this example we investigate the behavior of
our methodology when the class of convex sets used to fit the model differs
from that used to generate the true graph.  We consider three possibilities:
$\cA=\Alpha$ in $\rr^2$, $\cA=\Alpha$ in $\rr^3$ and $\cA=\Cech$ in $\rr^2$.
We performed two experiments: one when the graph is feasible for each the
three classes, and another example where the graph could be generated by only
two of the classes.

\begin{figure}[t!]
\begin{center}
  \includegraphics[height=100mm]{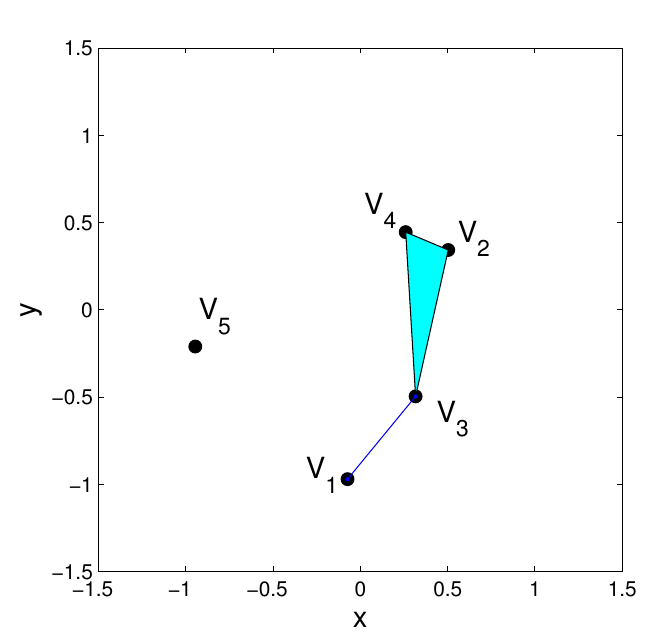}
  \caption{Graph encoding the Markov structure of the model given in
     \Eqn{Eq:ModelTalk1}.}\label{fig:ModelTalk1}
\end{center}
\end{figure}

First consider a model with junction tree factorization:
\begin{equation}\label{Eq:ModelTalk1}
  f_\theta(\bx)=\frac{\psi_\theta(x_1,x_3)\psi_\theta(x_2,x_3,x_4)\psi_\theta(x_5)}
                    {\psi_\theta(x_3)},
\end{equation}
whose conditional independence structure given by 
% a junction tree of the $1$-skeleton of the complex displayed in 
the graph of \Fig {fig:ModelTalk1}.  Again, the clique marginals are
specified as a Clayton copula with %association parameter
$\theta=4$.  We simulated $300$ samples from this distribution.

We fitted the model with each of the three classes of convex sets using the
Metropolis Hastings algorithm of \Sec{ss:samp} with random walk proposals on
$\BB \Dim$ (where $\Dim=2$ or $3$, depending on $\cA$).  \Alg {filtrationalg}
was used to enforce decomposability, using $r=0.40$ and $\eta=0.020$.  The
same exponential prior and uniform reflecting random-walk proposals for
$\theta$ were used as in Example \ref{Sec:ExampleIlust}.  Results of $1\,000$
samples after a burn-in period of $50\,000$ draws are summarized in
\Tab{Tab:ModelTalk1}.  Not surprisingly, the posterior mode coincided with
the true model in all three cases.

\begin{table}[t!]
\begin{center}
\begin{tabular}{|l|L|C|}
\hline\Strut Nerve &\text{HPP Models}&\text{Posterior}\\
\hline\Strut Alpha in $\rr^2$ 
%%%% &$\red{[2,3,4][1,3][5]}$ &$0.972$ \\
%%%% &$[2,3,4][1,2,3][5]$&$0.016$ \\
%%%% & $[2,3,4][1,2,3][3,5]$&$0.009$ \\
%% > Alpha in R^2   [2,3,4][1,3]  [5]   0.964
%% >                [2,3,4][1,3]  [1,5] 0.012
%% >                [2,3,4][1,3,4][5]   0.012
& \CQ{1,3}\CQ{2,3,4}\CQ{5}       & 0.964   \\
& \cq{1,3}\cq{2,3,4}\cq{1,5}     & 0.012   \\
& \cq{1,3,4}\cq{2,3,4}\cq{5}     & 0.012   \\
\hline\Strut Alpha in $\rr^3$ 
%%%% &$\red{\cq{2,3,4}\cq{1,3}\cq{5}}$   &$0.490$\\ 
%%%% &$\cq{1,3}\cq{2,3}\cq{3,4}\cq{5}$ &$0.103$\\
%%%% &$\cq{1,3}\cq{2,3}\cq{2,4}\cq{5}$ &$0.068$\\
%% > Alpha in R^3   [2,3,4][1,3]  [5]   0.982
%% >                [2,3] [3,4][1,3][5] 0.011
%% >                [2,3,4][1][5]       0.003
& \CQ{1,3}\CQ{2,3,4}\CQ{5}       & 0.982 \\
& \cq{1,3}\cq{2,3}\cq{3,4}\cq{5} & 0.011 \\
& \cq{1}\cq{2,3,4}\cq{5}         & 0.003 \\
\hline\Strut \Cech\ in $\rr^2$ 
%%%% &$\red{\cq{2,3,4}\cq{1,3}\cq{5}}$  &$0.607$\\ 
%%%% &$\cq{1,2,3,4}\cq{5}$      &$0.098$\\
%%%% &$\cq{1,2,3}\cq{2,3,4}\cq{5}$ &$0.061$\\
%% > Cech in R^2    [2,3,4][1,3][5]     0.595
%% >                [1,2,3,4] [5]       0.179
%% >                [1,2,3,4,5]         0.168
& \CQ{1,3}\CQ{2,3,4}\CQ{5}   & 0.595\\
& \cq{1,2,3,4}\cq{5}         & 0.179\\
& \cq{1,2,3,4,5}             & 0.168\\
\hline
\end{tabular}
\caption{Models with highest posterior probability.  The table is divided
  according to the class of convex sets used when fitting the model.  The
  true model is shown in bold.\label{Tab:ModelTalk1}}
%has $\cq{2,3,4}$, $\cq{1,3}$ and $\cq{5}$ as cliques.
\end{center}
\end{table}

\begin{figure}[t!]
\begin{center}
  \includegraphics[height=80mm]{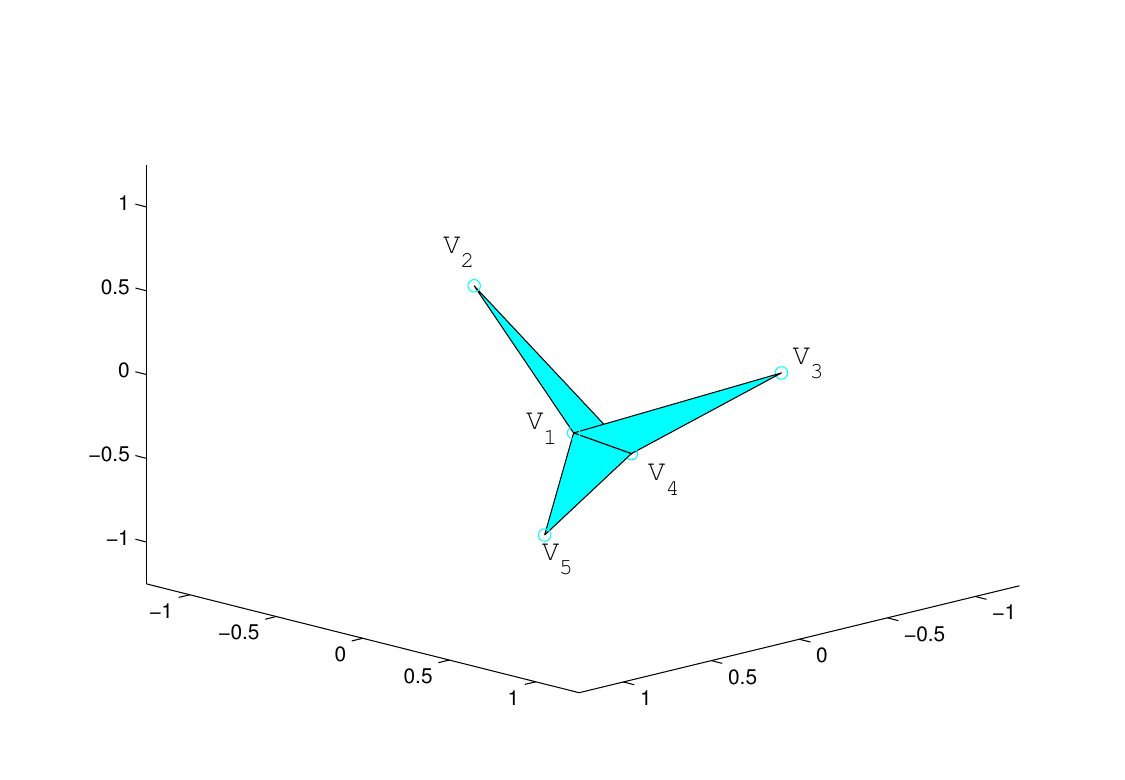}
  \caption{Graph encoding the Markov structure of the model given in
     \Eqn{Eq:ModelTalk2}.}\label{fig:ModelTalk2}
\end{center}
\end{figure}

The second model we considered has junction tree factorization:
\begin{equation}\label{Eq:ModelTalk2}
  f_\theta(\bx)=\frac{\psi_\theta(x_1,x_2,x_4)\psi_\theta(x_1,x_3,x_4)
                       \psi_\theta(x_1,x_4,x_5)}{\left(\psi_\theta(x_1,x_4)\right)^2}.
\end{equation}
The corresponding graph cannot be obtained from an Alpha complex in $\rr^2$,
but it is feasible for an Alpha complex in $\rr^3$ (\Fig {fig:ModelTalk2}) or
a \Cech\ complex in $\rr^2$.  Using the same Clayton clique marginals before,
we sampled $300$ observations from this distribution and fitted the model
using the three classes of convex sets.  Results from $1\,000$ samples after
a burn-in period of $75\,000$ are summarized in \Tab{Tab:ModelTalk2}.
Observe that for Alpha complexes in $\rr^2$, there is no clear posterior mode
(unlike the previous example, or Sections \ref {Sec:ExampleIlust} and \ref
{Sec:FacNerve}).  The posterior mode for the \Cech\ complex in $\rr^2$ and
the Alpha complex in $\rr^3$ both match the true model.

\begin{table}[t!]
\begin{center}
\begin{tabular}{|c|l|c|}
\hline\Strut 
Nerve &  HPP Models &  Posterior \\
\hline\Strut  Alpha in $\rr^2$
 %% %% & $\cq{1,2,3}\cq{1,3,4}\cq{1,4,5}$   &$0.383$  \\
 %% %% & $\cq{1,3,4}\cq{2,3,4}\cq{1,4,5}$   & $0.136$ \\
 %% %% & $\cq{1,2,3}\cq{1,2,4}\cq{1,4,5}$   & $0.104$  \\
 %% > Alpha in R^2   [1,2][1,3,4][1,4,5]   0.214
 %% >                [1,3,5][1,3,4][1,2,4]  0.115
 %% >                [3,4,5][1,3,4][1,2,4]  0.112
 & \cq{1,2}\cq{1,3,4}\cq{1,4,5}  & 0.214 \\
 & \cq{1,2,4}\cq{1,3,4}\cq{1,3,5}& 0.115 \\
 & \cq{1,2,4}\cq{1,3,4}\cq{3,4,5}& 0.112 \\
\hline\Strut  Alpha in $\rr^3$
 %% %% &$\cq{1,2,3}\cq{2,3,4}\cq{1,2,5}$ &$0.490$ \\ 
 %% %% &$\red{\cq{1,2,4}\cq{1,3,4}\cq{1,4,5}}$ &$0.103$ \\
 %% %% &$\cq{1,3,5}\cq{3,4,5}\cq{2,3,4}$ &$0.068$ \\
 %% > Alpha in R^3   [1,2,4][1,3,4][1,4,5] 0.976
 %% >                [1,4,5][1,2,3,4]      0.016
 %% >                [1,4,5][1,3][1,2,4]   0.009
 & \CQ{1,2,4}\CQ{1,3,4}\CQ{1,4,5} & 0.976 \\
 & \cq{1,2,3,4}\cq{1,4,5}         & 0.016 \\
 & \cq{1,2,4}\cq{1,3}\cq{1,4,5}   & 0.009 \\
\hline\Strut  \Cech\ in $\rr^2$
 %% %% &$\red{\cq{1,2,4}\cq{1,3,4}\cq{1,4,5}}$ &$0.607$\\ 
 %% %% &$\cq{1,2,4}\cq{1,3,4}\cq{1,3,5}$    &$0.098$\\
 %% %% &$\cq{1,2,4}\cq{1,3,4}\cq{4,5}$ &$0.061$\\
 %% > Cech in R^2    [1,2,4][1,3,4][1,4,5] 0.758
 %% >                [1,2,4][1,3,4][1,3,5] 0.177
 %% >                [1,2,3,4][1,4,5]      0.148
 & \CQ{1,2,4}\CQ{1,3,4}\CQ{1,4,5} & 0.758 \\
 & \cq{1,2,4}\cq{1,3,4}\cq{1,3,5} & 0.177 \\
 & \cq{1,2,4}\cq{1,3,4}\cq{4,5}   & 0.148 \\
\hline
\end{tabular}
\caption{Models with highest posterior probability, for each class of convex
  sets.  The true model (shown in bold) is unattainable for Alpha complexes
  in $\rr^2$. \label{Tab:ModelTalk2}}
%has $\cq{1,2,4}$, $\cq{1,3,4}$ and $\cq{1,4,5}$ as cliques.
\end{center}
\end{table}

\subsection{Comparative performance analysis with state-of-the-art}

We compared the performance of our method to Feature Inclusion Stochastic Search (FINCS), proposed by \cite{Scot:Carv:2008}, and the adaptive LASSO as described in the paper by \cite{Fan:Feng:Wu:2009}. The criterium for the comparison was given by estimating the counts of specific subgraphs of a set of graphs of size $50$. This is, we generated $100$ graphs of size $50$, and we computed the absolute counts for the following subgraphs: triangles, 4-cycles, 5-cycles, 3-stars, 4-stars, and 5-stars, for the true graph and for the estimated network by our method and its competitors, then we computed the absolute difference between the counts for the true and estimated graph and then divided by the count for the true graph. We denote by an $*$ the counts of induced subgraphs; this measures the performance of the methods under decomposability. We generated the set of true networks from an \ER\ random graph model with $\alpha=0.05$, so we are in a regime that can be called sparse.

To fit our method, we assumed an uniform distribution on the unit ball in $\rr^3$ for the vertex set and for the radius of each vertex an $\Ex(0.1)$ as priors. For the positions of the vertices we used the same proposals as in Examples \ref{Sec:ExampleIlust} - \ref{Sec:3Filtrations}. For each radius we implemented a mixture of random walks reflecting at $0$ as proposal. We used the nerve corresponding to the Cech complex. We set up FINCS with a probability of $0.1$ for resampling moves and a probability of $0.1$ for global moves. For the adaptive LASSO we set up $\gamma=0.5$ and $\Omega$ was initialized as the inverse of the sampled covariance matrix. Subgraphs were counted using the igraph command graph.subisomorphic.lad. When computing the normalized error, we adopted the convention $0/0=0$. For the induced subgraphs, we only compute the error of the Bayesian procedure and Lasso over the set of non-decomposable graphs. For the simulation displayed in Table \ref{tab:Comparing} all $100$ graphs were non-decomposable.

\begin{table}[b!]
\begin{center}
\begin{tabular}{|L C C C|}
\hline\Strut
\text{Subgraph}&\text{Bayes}&\text{FINCS}&\text{Lasso} \\
\hline\Strut

\text{triangles} &     0.083   \pm 0.07          &     0.125  \pm  0.17        &     0.208    \pm 0.12          \\
\text{4- cycles}  &    0.062  \pm 0.10          &      0.123  \pm  0.16        &    0.166    \pm 0.10           \\
\text{5- cycles} &     0.086  \pm 0.07          &     0.112    \pm  0.14       &   0.124     \pm 0.08          \\
\text{3- stars}  &      0.103   \pm 0.08        &      0.139    \pm  0.12       &   0.211     \pm 0.08          \\ 
\text{4- stars} &       0.087   \pm 0.08       &      0.096     \pm  0.16    &    0.115        \pm 0.11      \\
\text{5- stars}  &      0.201   \pm0.10       &     0.183      \pm   0.14       &  0.174       \pm 0.06          \\ 
\text{4- cycles}  * &  0.146  \pm 0.12      &      0.930  \pm  0.07          &    0.229    \pm 0.12           \\
\text{5- cycles}  * &  0.128  \pm 0.13      &     1      \pm  0                   &   0.174     \pm 0.08          \\

\hline
\multicolumn{2}{c}{}
\end{tabular}
\caption{Estimated normalized errors for counts of specific subgraphs for our method, FINCS and adaptive LASSO.}\label{tab:Comparing}
\end{center}
\end{table}

Results are summarized in Table \ref{tab:Comparing}. Our method incurred into less errors on average compared to our competitors for almost all subgraphs. The exception was the $5-$star. We also observed that FINCS outperformed the LASSO for almost all regimes, with exception of the $5-$star.

\begin{figure}[t!]
\begin{center}
  \includegraphics[height=100mm,trim=120 0 100 0]{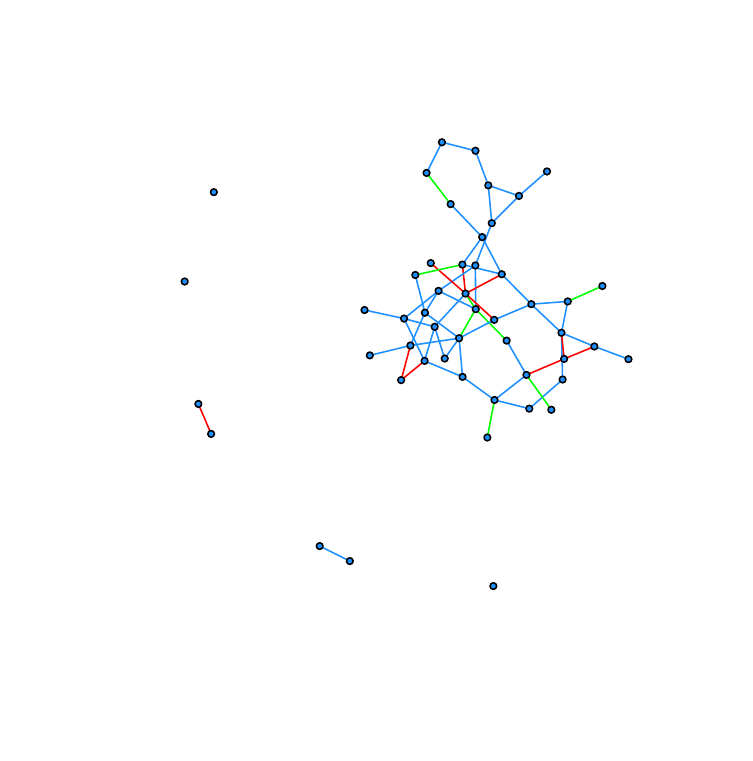}
  \caption{Here we compare the true model, which was sampled from the ERGM used in Table \ref{tab:CompareStruc}, and the fitted model, using our method (again, as in Table \ref{tab:CompareStruc}).  The edges that were added by our method (with respect to the true graph) are highlighted in red. The edges that were deleted by our method are highlighted in green.}\label{fig:SimulERGM}
\end{center}
\end{figure}

We performed another simulation, now assuming an exponential random graph model for the true graph. The simulation of the true graphs was implemented using the R package statnet. We used the formula  edges+triangles+kstar(3), with coef=c(-3.2,0.95,0.005); this specification encourages the presence of triangles and 3-stars. These choices produce graphs with twice the 3-stars and 3 times more triangles than an  \ER\  with $\alpha=0.05$, while having approximately the same density. The objective of this experiment is to investigate the behavior of our method when the true graph has more structure than the typical realization of an \ER\  model. Results are summarized in Table \ref{tab:CompareStruc}; they are similar to what was observed in the previous experiment. We observed an improvement regarding the counts of triangles, this is not surprising since geometric random graphs tend to have more triangles than realizations from an \ER\ model. In Figure \ref{fig:SimulERGM} we compare the true model (as generated from the ERGM just described) and the fitted model for a single realization. In this regime, graphs tend to be non-decomposable. We estimated the proportion of decomposable graphs from a sample of 1,000 networks sampled from the ERGM used to obtain the simulation in Table \ref{tab:CompareStruc}, and obtained 0.302 as the result. For the simulation displayed in Table \ref{tab:CompareStruc}, $72$ out of the $100$ graphs were non-decomposable.

\begin{table}[t!]
\begin{center}
\begin{tabular}{|L C C C|}
\hline\Strut
\text{Subgraph}&\text{Bayes}&\text{FINCS}&\text{Lasso} \\
\hline\Strut

\text{triangles} &     0.071   \pm 0.04          &     0.134   \pm 0.14         &     0.217   \pm 0.09          \\
\text{4- cycles}  &    0.067   \pm 0.06         &      0.121  \pm 0.12          &    0.154   \pm  0.04          \\
\text{5- cycles} &     0.075   \pm 0.09          &     0.118  \pm 0.11           &   0.131    \pm 0.13          \\
\text{3- stars}  &      0.092   \pm 0.07         &      0.144  \pm 0.14           &   0.236    \pm 0.12          \\ 
\text{4- stars} &       0.086   \pm 0.09         &      0.115  \pm 0.15           &    0.117    \pm 0.10         \\
\text{5- stars}  &      0.214   \pm 0.09         &     0.122   \pm 0.13           &  0.121      \pm 0.06          \\ 
\text{4- cycles}  * &  0.152  \pm 0.12      &        0.720  \pm  0.07          &    0.214    \pm 0.10           \\
\text{5- cycles}  * &  0.133  \pm 0.10      &       0.720   \pm  0.13          &   0.163     \pm 0.08          \\

\hline
\multicolumn{2}{c}{}
\end{tabular}
\caption{Estimated normalized errors for counts of specific subgraphs for our method, FINCS and adaptive LASSO. The true graphs were sampled from an ERGM that encouraged the presence of triangles and 3-stars.}\label{tab:CompareStruc}
\end{center}
\end{table}

\subsubsection{Scalability}

Here we discuss scalability of the proposed method and of the competing methods ( \cite{Fan:Feng:Wu:2009}, \cite{Scot:Carv:2008}) to better appreciate the cost incurred in producing the errors in Tables \ref{tab:Comparing} and \ref{tab:CompareStruc}. Since the implementation of the proposed and competing methods available to us are in different programming languages, which influence greatly the actual runtime, we outline such a discussion in theory, in terms of key quantities that influence scalability, including number of nodes $n$, number of edges $m$, and number of cliques $k$. We also distinguish two tasks: the task of estimating the parameter of a model-based representation of a (hyper)graph, and the task of generating $b$ (hyper)graphs from an estimated model-based representation. 
%For reference, also note that the literature is typically concerned with a sparse graph regime, where $m=O(n)$, and a dense graph regime, where $m=O(n^2)$.

Regarding the task of estimating parameters from an observed (hyper)graph, the proposed methods requires estimating parameters $\set{V},r$, and the estimation complexity scales as $O(S(n)+k^3)$, where $S(n)$ denotes the complexity of matrix multiplication.\footnote{The parameter $k$ in this case is actually the number of prime components, but this quantity is typically in the same order of magnitude as the number of cliques.} The method by  \cite{Scot:Carv:2008} requires estimating parameters $\cG$, and the estimation complexity scales as $O(S(n)+k)$. The lasso requires estimating parameters $\Sigma^{-1}$, and the estimation complexity scales as $O(n^3)$. 

Regarding the task of generating $b$ (hyper)graphs from an estimated model-based representation, the proposed methods scales as $O(bn^2)$. 
This is because the complexity of computing the 2-skeleton of the Cech complex scales as $O(n^2)$. Alternative representations lead to different scaling: computing the Delaunay triangulation scales as $O(n\log n)$, and computing the the Alpha complex scales as $O(n^2)$.
The method by \cite{Scot:Carv:2008} scales as $O(bnm)$, where typically $m$ is much larger than $n$. For the lasso, this statement does not apply.

To illustrate how our method scales up with respect to the number of variables, we ran the experiment summarized in Table \ref{tab:Scale}. We obtained the graphs from the nerves of \Cech\ complexes and employed the method proposed by \citep {Atay:Mass:2005} to compute the normalizing constants of non-decomposable graphs. The MCMC was performed on a 2.5GHz desk computer with 4GB of RAM. Our method was implemented using Matlab (MathWorks).
\begin{table}[ht]
\begin{center}
\begin{tabular}{|L C C R|}
\hline\Strut
\text{N Variables}&\text{Burn-in}&\text{N Iterations}&\text{Time} \\
\hline\Strut

4    &     50,000          &     10,000         &       7m \\ %438.231 Sec.          \\
40  &     50,000         &      10,000         &    2h~11m \\ %7,869.352 Sec.         \\
400 &    50,000         &      10,000         &  3d~2h~17m \\ % 27,385.264 Sec.         \\
\hline
\multicolumn{2}{c}{}
\end{tabular}
\caption{This experiment illustrates how our method scales up with respect to the number of variables.}\label{tab:Scale}
\end{center}
\end{table}

\subsection{Real data analysis}

\subsubsection{Fisher's Iris data}

We applied our method to Fisher's Iris data set. Variables include: sepal length (1), sepal width (2), petal length (3), and petal width (4). The objective is to find the conditional independence structure given a family of distributions for the likelihood (e.g. Gaussian, Clayton copula). A summary of the distribution of this data set is given in Table \ref{tab:Iris}

\begin{table}[b!]
\begin{center}
\begin{tabular}{|L C C C C|}
\hline\Strut
\text{Variable}&\text{SL}&\text{SW}&\text{PL}&\text{PW}\\
\hline\Strut

\text{Sepal length} &  0.4043    &     &      &     \\
\text{Sepal width}  &   0.0938   & 0.1040    &      &     \\
\text{Petal length} &  0.3033     & 0.0714     &  0.3046     &     \\
\text{Petal width}  &  0.0491     & 0.0476     &  0.0488     & 0.0754   \\ 

\hline
\multicolumn{2}{c}{}
\end{tabular}
\caption{Estimated variances and covariances for the Iris data.}\label{tab:Iris}
\end{center}
\end{table}

We first describe the specification we used for the random geometric graph, then we will make our choices for the Hyper-Markov law and the likelihood explicit. For the RGG we assumed an uniform distribution of the vertices on the unit ball in $\rr^3$ and for the radius of each vertex an $\Ex(0.1)$, distribution was assumed. For the positions of the vertices we used the same proposals as in Examples \ref{Sec:ExampleIlust} - \ref{Sec:3Filtrations}. For each radius we used a mixture of random walks reflecting at $0$ as proposal. For the likelihood function and Hyper-Markov law we used the following specifications:  A multivariate normal distribution for the likelihood with an HIW as Hyper-Markov law. 
%For the second one we used a factorization of Clayton copulas (as in Examples \ref{Sec:ExampleIlust},\ref{Sec:FacNerve}, and \ref{Sec:3Filtrations}) with an exponential distribution for $\theta$, the copula parameter,  as Hyper-Markov law.

For our choice for the likelihood and Hyper-Markov law, we adopted the same values for the hyperparameters as  \citep{Rove:2002} did, this is, the prior for the precision matrix centered at  $I$ and $\delta=3$. We used the method proposed by \citep {Atay:Mass:2005} to deal with the normalizing constants of non-decomposable graphs.  We ran the MCMC with $300\,000$ iterations as burn-in and kept the last $10\,000$ for analysis.

Results for the first choice are summarized in Table \ref{tab:FactRes1}. Here we display the $9$ models with highest posterior probability. All the posterior probability was concentrated in this models. Our posterior mode coincides with the one reported by \citep{Rove:2002}, but we obtained different results for the rest of the models. We attribute this difference to the fact that we used different priors for graph space; ours being non-uniform.

\begin{table}[b!]
\begin{center}
\begin{tabular}{|L C|}
\hline\Strut
\text{Maximal Simplices}&\text{Posterior Probability}\\
\hline\Strut
\cs{1,2}\cs{1,3}\cs{2,4}\cs{3,4}      & 0.3465\\
\cs{2}\cs{1,3}\cs{3,4}       & 0.2835\\
\cs{1,3}\cs{2,3}\cs{4}      & 0.1999\\
\cs{1,2}\cs{1,3}\cs{4}      & 0.1540\\
\cs{1,2}\cs{1,4}\cs{2}      & 0.0116\\
\cs{1,4}\cs{2}\cs{3,4}      & 0.0026\\
\cs{1,2}\cs{1,3}\cs{3,4}      & 0.0016\\
\cs{1,3}\cs{2,3}\cs{3,4}      & 0.0002\\
\cs{1,4}\cs{1,3}\cs{2,3}\cs{3,4}      & 0.0001\\
\hline
\multicolumn{2}{c}{}
\end{tabular}
\caption{Highest posterior factorizations with uniform prior and Gaussian distribution for Fisher's Iris data set.}\label{tab:FactRes1}
\end{center}
\end{table}

%\begin{table}[ht]
%\begin{center}
%\begin{tabular}{|L C|}
%\hline\Strut
%\text{Maximal Simplices}&\text{Posterior Probability}\\
%\hline\Strut
%\cs{2,3,4}\cs{1}     & 0.3673\\
% \cs{2,3,4}\cs{1,4}\cs{1,3}       & 0.3227\\
% \cs{2,3,4}\cs{1,2}\cs{1,3}       & 0.1726\\
%  \cs{2,3,4}\cs{1,2}\cs{1,4}       & 0.1374\\
%\hline
%\multicolumn{2}{c}{}
%\end{tabular}
%\caption{Highest posterior factorizations with uniform prior and Clayton copula for Fisher's Iris data set.}\label{tab:FactRes2}
%\end{center}
%\end{table}

%Results for the second choice of likelihood and Hyper-Markov law are summarized in Table \ref{tab:FactRes2}. Here all the posterior probability is concentrated in $4$ models. Since the family of distributions used to describe the data is different from the previous one, it was expected to obtain different results. It is worth mentioning that the hyperedge $[2,3,4]$ is a factor in all four models.

\begin{table}[b!]
\begin{center}
\begin{tabular}{|L C|}
\hline\Strut
\text{Maximal Simplices}&\text{Frequency as posterior mode}\\
\hline\Strut
\cs{1,2}\cs{1,3}\cs{2,4}\cs{3,4}      & 0.68\\
\cs{1,3}\cs{2,3}\cs{4}      & 0.16\\
\cs{1,2}\cs{2,3}\cs{2,4}      & 0.08\\
\cs{1,4}\cs{2,4}\cs{3,4}      & 0.04\\
\cs{1,4}\cs{1,3}\cs{2}        & 0.04\\
\hline
\multicolumn{2}{c}{}
\end{tabular}
\caption{Results from missing data simulation: We fit our model the same way as in Table \ref{tab:FactRes1}, but leaving one row out each time and imputing it as missing data. In this table we report the frequency by which each factorization appeared as posterior mode.}\label{tab:FactRes3}
\end{center}
\end{table}

We conducted another simulation, where we assessed the robustness of the inference for the Gaussian graphical model via tools from missing data. We fit our model deleting one row of the data at a time (this is, we fit our model 50 times over incomplete data sets) and imputed the missing data using the conditional distribution of the observed data. Results are summarized in \ref{tab:FactRes3}. For each of this simulations we computed the average distance between the predicted vector and the observed one. For FINCS, the average distance between the predicted and observed vectors (across the 50 simulations) was $3.69 \pm 1.28$, while for our method it was $3.22 \pm 1.15$. This is not surprising, since, in contrast to FINCS, we consider non-decomposable models. 

\subsubsection{Daily exchange rates data}

Following \cite{Hern:Lloy:Hern:2013}, we considered daily exchange rates of eight currencies (Swiss franc, Australian dollar, Canadian dollar, Norwegian krone, Swedish krona, Euro, New Zealand dollar and British pound) with respect to the U.S. dollar. The data set consists of 717 observations, from 1 Dec., 2011, to 29 Aug., 2014. Clearly these observations are not \emph{iid}, but we will not take this into account in the modeling. What makes this application interesting is the presence of a non-trivial missing data pattern. The data was downloaded from yahoo.finance.com via the R library tseries.

\begin{figure}[t!]
\begin{center}
  \includegraphics[height=100mm]{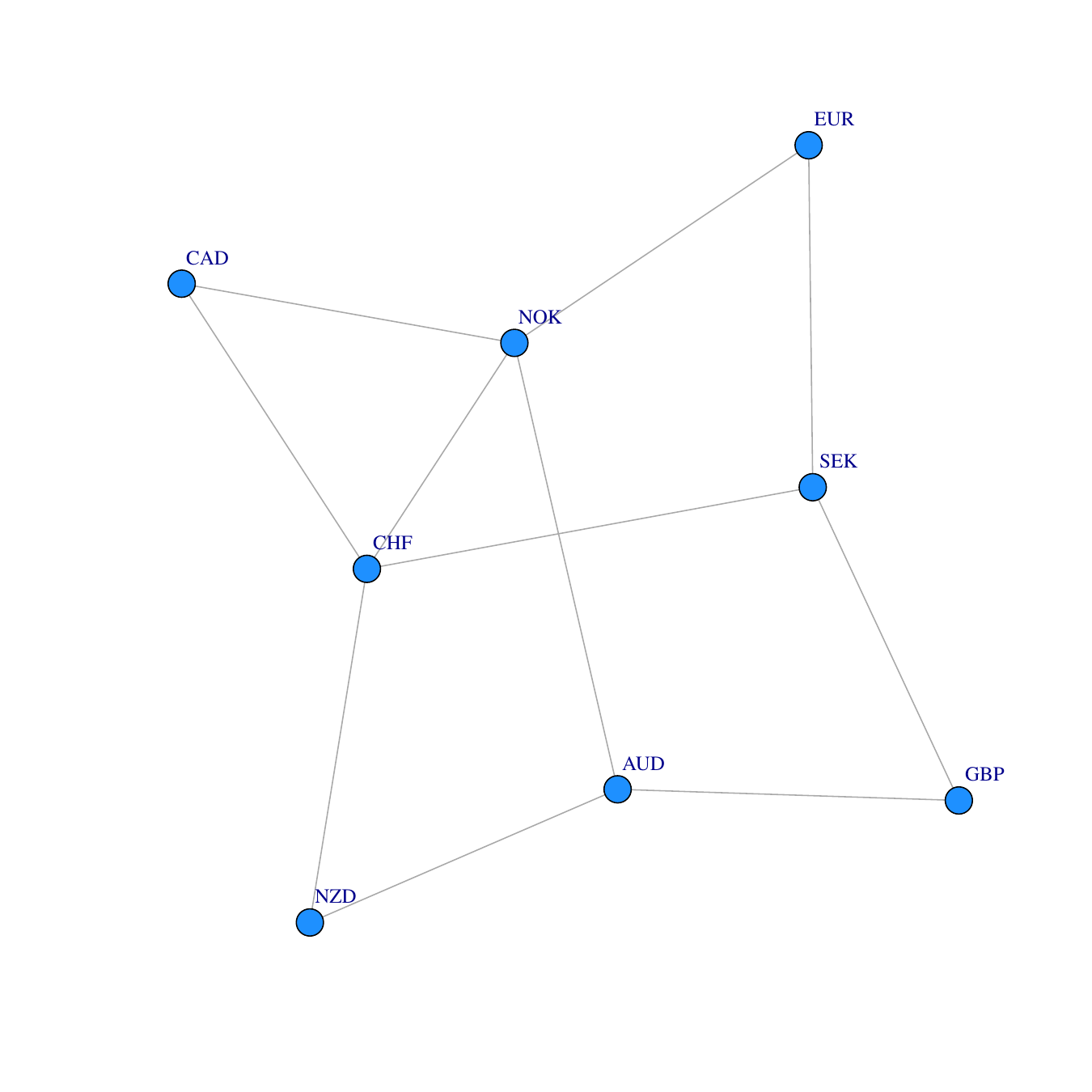}
  \caption{Posterior mode of our method when applied to daily exchange rates (with respect to US Dollar) from 1 Dec., 2011 to 29 Aug., 2014. Here CHF denotes the Swiss franc, AUD, the Australian dollar, CAD, the Canadian dollar, NOK, the Norwegian krone, SEK, the Swedish krona, EUR, the Euro, NZD, the New Zealand dollar, and GBP, the British pound. }\label{fig:Finance}
\end{center}
\end{figure}

We applied our method to these data. We assumed a uniform distribution for the vertices in the unit ball in $\rr^3$, and that the nerve was computed from the intersection pattern of balls of different sizes. We assumed a HIW as the hyper-Markov law and a multivariate normal as the likelihood. We kept the simulated values from 5,000 iterations after 25,000 iterations of burn-in. Missing data were assumed missing at random and imputed from the model. The posterior mode is illustrated in Figure \ref{fig:Finance}; it has 0.87 posterior probability. This model is non-decomposable.

\section{Discussion}
% Section 6

% EDO : THIS DISCUSSION IS GOOD. WE SHOULD MIGRATE SOME OF THE FUTURE RESEARCH IDEAS TO THE CONCLUDING REMARKS. ALSO, SOME OF THESE DISCUSSION POINTS WILL BE REPLACED BY MORE RECENT DEVELOPMENT - EG, SEE SPACE OF GRAPHS. I AM ADDING THE POINTS WE DISCUSSED BELOW, WE CAN THINK ABOUT HOW TO ORGANIZE THEM AT A LATER POINT, ONCE YOU HAVE TEXT FOR THEM

%New points to be fleshed out.

%1. 
  
%2. 

%4. 

%5. \\

%{*** old discussion starts below. keep it, but make sure to purge the points we are flashing out in the new discussion ***}

\bigskip
In this article we present a new parametrization of graphs by associating
them with finite sets of points in $\rr^\Dim$.  This perspective supports the
design of informative prior distributions on graphs using familiar
probability distributions of point sets in $\rr^\Dim$.  It also induces new
and useful Metropolis\slash Hastings proposal distributions in graph space
that include both local and global moves.  As suggested by Helly's Theorem
\citep {Edel:Hare:2009} characterizing the sparsity of intersections of
convex sets in $\rr^\Dim$, this methodology is particularly well suited for
sparse graphs.  The simple strategies presented here generalize easily to
more detailed and subtle models for both priors and M/H proposals.

Our construction leads to MCMC that naturally instantiate local and global moves in graph (and hypergraph) space. The proposals that produce small perturbations on the vertex set will produce local moves with high probability, while proposals that consist in resampling a subset of the vertex set will produce, with high probability, a global move (See Figure \ref{Fig:Global}). 

An interesting feature of our approach is that the distribution on the space
of graphs is modeled directly before the application of any specific hyper
Markov law, in contrast to standard approaches in which it is the hyper
Markov law that is used to encourage sparsity or other features on the graph.
We think that working with the space of graphs explicitly opens a lot of
possibilities for prior specification in graphical models, therefore, it is a
perspective worth further study.

While coupling the focus on the first two moments and a graph representation of pairwise dependencies among variables are not restrictive modeling choices in the Gaussian graphical model framework, they become restrictive when working with graphical models in general. However, likely because of convenience, these assumptions are seldom challenged in the graphical models literature. Here we develop a geometric construction where dependence relations of higher orders can be conveniently encoded within a hypergraph. For a state of the art treatment of parametrizations of decomposable graphs, see \citep{Cron:West:2012}.

Connected to the point above, while decomposability plays an important role in graphical models in general, it does not play any role in our construction. This is because decomposability is relevant for computations at the Hyper-Markov law level, while our construction and results are at the level of the prior on graph space.

About the space of graph where our construction puts positive mass. At this point we have two results and a conjecture for characterizing the space of feasible graphs (hypergraphs) when we consider the nerve of a set of balls in $\rr^3$ with different radii (This is the case that should have the largest support for the distribution of graphs.)
 
Theorems: 
    (i) Any graph can be embedded in $\mathbb{R}^3$. This is a well known result. One method that computes such embedding is the Book Algorithm, proposed by \citep{Kainen:1974fk} and \citep{Ollm:1973}.
   (ii) Any graph can be linearly embedded in $\mathbb{R}^3$. By `linearly embedded' it is meant that the segments joining the vertices are straight lines. This is a particular case of a more general result, that every k-dimensional simplicial complex can be geometrically realized in $2k+1$ dimensions.  In this case case, $k=1$. See Section 3.1 of  \citep {Edel:Hare:2008}.

 Conjecture: 
  (iii) Such a linear embedding can be achieved using a random geometric graphs construction using balls with different radii.

Interesting questions and extensions of this idea include: (1)~achieving a
deeper and more detailed understanding of the subspace of graphs spanned by
different specific filtrations; (2)~designing priors to control the
distributions of specific features of graphs such as clique size or tree
width; (3)~modeling directed acyclic graphs (DAGs), and (4)~concrete
implementation of novel Markov structures based on nerves.

%The proposed methodology requires that the graph is computed from a
%filtration which we obtain from nerves.  Ideally we would like
%to characterize the model space generated from the nerves in
%terms of topological concepts.  In the case of decomposable
%graphs we have not been able find a purely topological 
%characterization which makes the analysis of the model
%space spanned by a specific filtration difficult.  Although 
%we do have some insight about which graphs cannot be generated by 
%a specific filtrations more formal understanding of this 
%aspect of the problem is needed.

This methodology generates only graphs that are feasible for the particular
filtration chosen.  Although we do have some insight about which graphs can
and cannot be generated by a specific filtrations, a more complete and formal
understanding of this aspect of the problem would be useful.

We used very simple prior distributions for the purpose of illustrating the
core idea of the methodology.  It is natural in this approach to incorporate
tools from point processes into graphical models to define new classes of
priors for graph space.  Future developments in our research will involve a
range of repulsive and cluster processes.

%We used simple marginal models to illustrate our ideas.  
%One extension is to incorporate approaches using hyper Markov laws 
%to couple our models for the graph with richer models on the marginals.
%An example is to incorporate more sophisticated models
%for the clique marginals (such as Bernstein copulas).  

The parametrization we propose can be used to represent Markov structures on
DAGs, but the strategies for obtaining such graphs from nerves will be
different and will establish stronger connections between Graphical Models
and computational topology.

The present work is related to that of \citet {Pist:Wynn:etal:2009} in which
a nerve of convex subsets in $\rr^\Dim$ is used to obtain Markov structures
for a distribution, an extension of the abstract tube theory of \citet
{Naim:Wynn:1997}.  This new perspective allows for constructions that
generalize the idea of junction trees.  By modifying our methodology
according to this framework (personal communication from H.\ Wynn suggests
that this is feasible) we hope to fit models that factorize according to
those novel Markov structures.

Another possible extension of this work is to discretize the set from which
the vertex set is sampled (\eg, use a grid).  Such discretization may improve
the behaviour of the MCMC; it would also allow the use of a nonparametric
form for the prior on the vertex set, leading to more flexible priors on
graph space.

While we illustrated the new construction for Bayesian inference, in a situation where we observe high-dimensional vectors and we want to infer the dependencies among their components, the proposed construction can be easily used to build a likelihood in situations where we have direct observations about the facets of hypergraph. Such observations occur naturally in applications; just think of how one would encode the structure among individuals revealed by pictures on Facebook. Each picture has one or more people. A picture with three people is naturally encoded as a 3-facet, rather than as 3 individual edges, as currently done, arguably for lack of likelihood models for hyper graphs.

%\section*{Acknowledgments}
%We are grateful to Herbert Edelsbrunner, John Harer, and Henry Wynn for helpful conversations.
%%
%This work was partially supported by National Science Foundation grants DMS--0732260, DMS--0635449, DMS-0757549, and IIS-1149662, by National Institute of Health grants NIH R01 CA123175, R01 GM096193, and NIH P50-GM081883, by the Army Research Office grant MURI W911NF-11-1-0036, and by an Alfred P. Sloan Research Fellowship.  
%%
%Any opinions, conclusions or recommendations expressed in this material are those of the authors and do not necessarily reflect the views of the NSF, the NIH, or the ARO.

%EMA is partially supported by NSF CAREER award , ARO MURI award W911NF-11-1-0036, and an Alfred P. Sloan Research Fellowship

%\bibliography{statjour-full,comptop}
\section*{References}
\vspace{-35pt}
\bibliographystyle{plainnat}
\bibliography{comptop}

\appendix
\section{Filtrations and Decomposability in Random Geometric Graphs}\label{Sec:Filtration}
% EDO : THIS APPENDIX NEEDS A BETTER TITLE. ALSO, I DID NOT PUT ANY TEXT AT THE END OF SECTION \REF{s:rgg} TO POINT HERE. PLEASE ADD IT.

A \emph{filtration} is a sequence of simplicial complexes that properly
include their predecessors:
\begin{definition}[Filtration]
  A \emph{filtration} for a simplicial complex $\cK$ is a sequence
  $\sL=\set{\cK^0,\cK^1,\ldots , \cK^k}$ of simplicial complexes
  such that
\begin{displaymath}
  \varnothing=\cK^0 \subset \cK^1 \subset \ldots \subset \cK^k=\cK,
\end{displaymath}
with all inclusions proper.  
\end{definition}

Filtrations are commonly used in computational geometry and topology to
construct efficient algorithms for computing specific nerves including the
Alpha complex \citep {Edel:Muck:1994}.  The simplicial complexes constructed
in \Sec {s:rgg} from families of convex sets lead to
filtrations as the convex sets are enlarged, by increasing the radius
parameter $r$.
%
% Fix a finite set $\cV\subset\rr^\Dim$ of $\Num\in\nn$ points and generate
% the nerve for $\cV$ with radius $r$ for each $r\ge0$ and, from this, the
% associated \Cech\ or Alpha complexes.  Because complexes only grow with
% increasing $r$, this induces a
% filtration. %since simplices are never removed
%
% \uhoh{fix this} An incremental approach to constructing a simplicial
% complex from a vertex set $\cV$ and a radius $r$ is to compute the nerve
% for $\cV$ with radius $0\le s<r$ with $s$ increasing from $0$ to $r$.  This
% procedure adds sets to the complex as $s$, sets are never removed.  This is
% true since $\cap_{j=1}^\Dim {\cal A}(V_j,s_1)=\varnothing$ implies
% $\cap_{j=1}^\Dim {\cal A}(V_j,s_2)=\varnothing$ for all $s_2>s_1\geq 0$.  A
% filtration is used as a structure to keep track of the addition of sets as
% the radius increases.  In computational geometry and topology filtrations
% are used to construct efficient algorithms for computing specific nerves
% such as the Alpha complex (this follows because $\cap_{j=1}^\Dim
% \cA(V_j,s)=\varnothing$ implies $\cap_{j=1}^\Dim \cA(V_j,r)=\varnothing$
% for all $0\le r < s$).  (this follows because the convex sets $B_{v,r}$ and
% $B_{v,r}\cap C_V$ only grow with increasing $r$).
%
% Filtrations serve as convenient structures to keep track of the addition of
% sets as the radius increases.  We will use filtrations to relate the
% concept of simplicial complexes to weak decomposability.  Note that the use
% of nerves to represent graphical models does not require that the graphs
% have to be weakly decomposable.

\begin{figure}[htb]
\begin{center}
  \includegraphics[width=\textwidth]{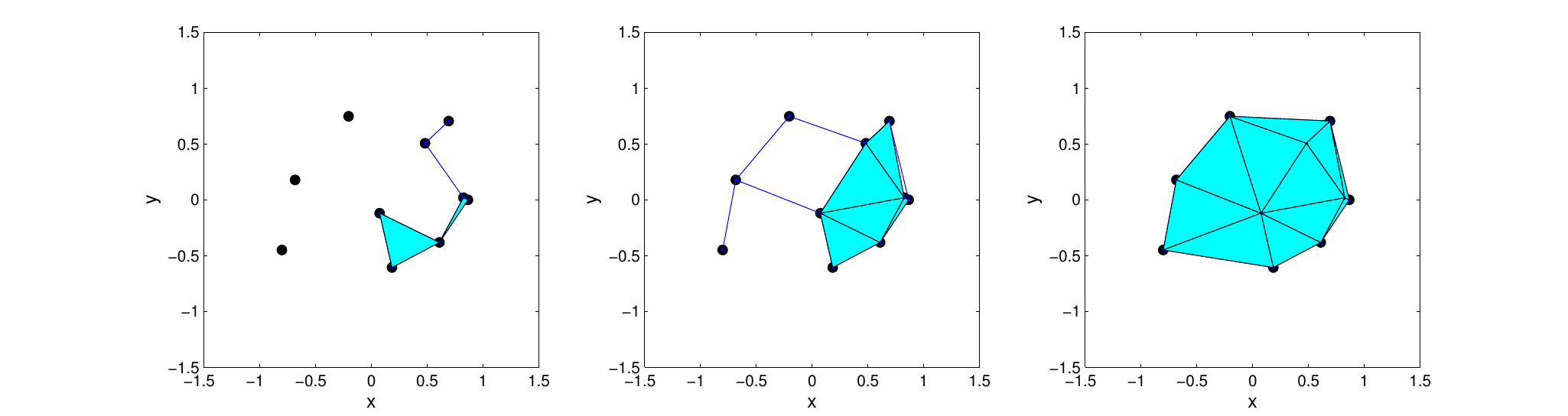}
  \caption{Filtration of Alpha complexes, (a) $r=\sqrt{0.10}$, (b)
    $r=\sqrt{0.20}$ and (c) $r=\sqrt{0.75}$.
    \label{fig:Filtration}}
\end{center}
\end{figure}

Although much of the graphical models literature focuses on Markov structures
derived from decomposable graphs, those constructed in \Sec
{s:rgg} from $1$-skeletons of \Cech\ and Alpha complexes need
not be decomposable (see \Fig {fig:Filtration}).  In \Alg
{filtrationalg} we present an adaptation of this construction that generates
 decomposable graphs, for use in applications that require them.  In
Sections \ref {sec:StatModel} and \ref{Sec:Simul} we present methodology and
examples for both decomposable and unrestricted model spaces.

The central idea for generating decomposable graphs from filtrations is to
note that the complex $\cA(\cV,0)$ for $r=0$ is disconnected and hence 
decomposable; as the radius $r$ increases, if one adds edges \emph{only} if
the resulting graph is decomposable, then decomposability will hold
for all $r\ge0$.  This procedure is formalized in \Alg{filtrationalg}.

\begin{algorithm}\label{filtrationalg}
\SetKwInOut{Input}{input} \SetKwInOut{Return}{return}

\vspace{1mm}
\Input{\mbox{a filtration
     $\sL=\set{\cK^0,\cK^1,\ldots , \cK^k}$}} \vspace{1mm}
\Return{\mbox{a weakly decomposable graph $\cG$}} \vspace{1mm}

$m=0$;\quad $i=0$;\quad $\cG_0 = \varnothing$\;

\While{$i < k \mathbf{~and~} \cK^{i+1} \neq \varnothing$}{
   $\tau_i =\set{ \kappa \in \cK^{i+1} \setminus \cK^i
   \mbox{ such that }  |\kappa|=2 }$;
   \tcp{the edges in the set difference and denote $\tau_{i,s}$
        as  the $s$th edge}\
   \If{$\tau_i \neq \varnothing$}{
      $P=|\tau_i|$\;
      \For{$s=1$ \KwTo $P$}{
         $\cG' = \cG_m \bigcup \tau_{i,s}$ ;  \tcp{propose adding the edge}
         \If{$\sharp(\cG')  < \sharp(\cG_m)$
             \tcp{fewer connected components?}} {
             $\cG_{m+1} = \cG'$ ;  \tcp{accept the proposal}\
             $m=m+1$\;}
         \Else{
            $\cq{c_i} = \sC(\cG_m)$ ; \tcp{the cliques }
            $\cq{s_i} = \sS(\cG_m)$ ; \tcp{the separators}
            $[v_1,v_2] = \tau_{i,s}\quad$; \tcp{the proposed edge}\
            \If{ $\exists i \neq j,k$ with $v_1 \in c_i, v_2 \in c_j$
                 and $c_i \bigcap c_j = s_k$ }
                 { $\cG_{m+1} = \cG'$ ; \tcp{accept the proposal}\
                   $m=m+1$\;
                 }
         }
      }
   }
   $i=i+1$\;
}
$\cG = \cG_m$

\caption{The algorithm takes as input a filtration of $k$ abstract complexes
  and returns a decomposable graph that is a subset of the $1$-skeleton of
  the $k$-th complex.}
\end{algorithm}

\begin{prop}
  The graph $\cG$ produced by \Alg{filtrationalg} is decomposable.
\end{prop} 

\begin{proof} The algorithm is initialized with the decomposable empty graph,
and decomposability is tested with each proposed addition of an edge
(\ie, a $1$-simplex in $\sL$).  The decomposability test is taken from
%\citeauthor {Giud:Gree:1999} \cite [Theorem $2$] {Giud:Gree:1999}.  
\citep [Theorem $2$] {Giud:Gree:1999}.  Since only finitely many edges may
possibly be added, $\cG$ is decomposable by construction.
\end{proof}

\begin{figure}[htb]
\begin{center}
  \includegraphics[height=70mm]{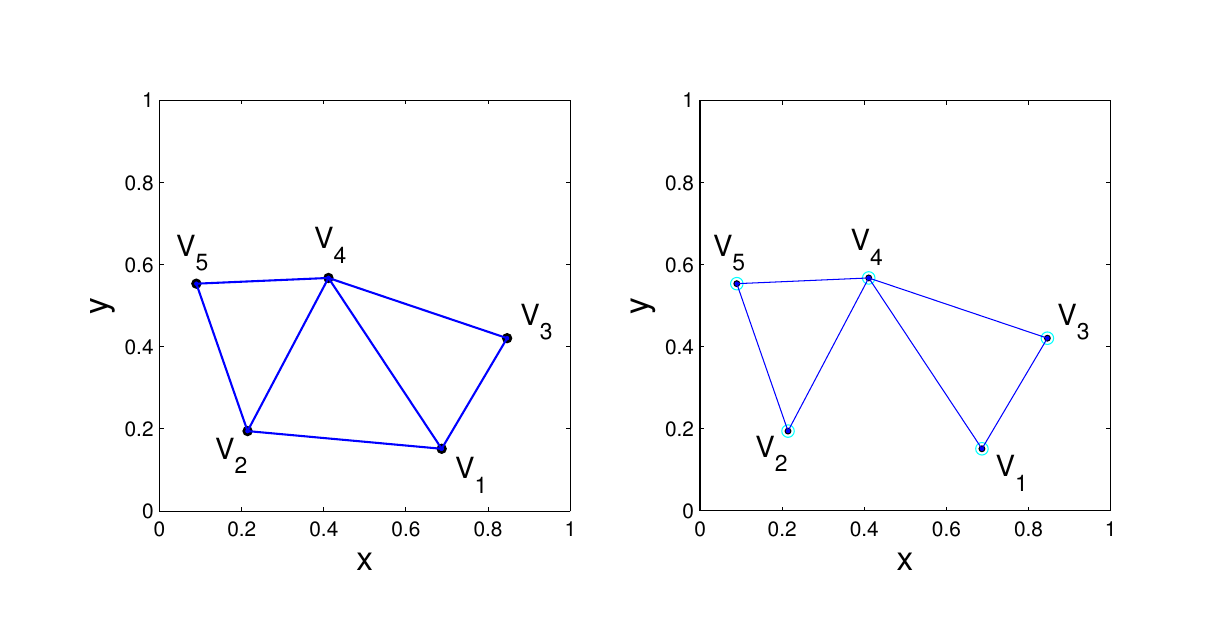}
  \caption{(a) Proximity graph computed from the vertex set given 
    in \Tab{tab:Alg1_Vertex}. (b) The decomposable graph computed 
    from the same vertex set using Algorithm \ref{filtrationalg}. 
    The edge $(1,2)$ is not included in the decomposable graph.}\label{fig:Alg1_Example}
\end{center}
\end{figure}

\subsection{Example}

We first illustrate the algorithm on a simple example, based on the five
points in $\rr^2$ shown in \Fig{fig:Alg1_Example} and given in
\Tab {tab:Alg1_Vertex}. The graph induced by a
\Cech\ complex with $r=0.5$ will not be decomposable.  \Tab
{tab:Alg1_Example} presents the evolution of cliques and separators with
increasing $r$, as edges are proposed for inclusion in \Alg{filtrationalg}.
The first three proposed edge additions are accepted, but the proposal to add
edge $(1,2)$ at radius $r=0.474$ is rejected, since
% in the filtration are added by $r=.379$ since this reduces the number of
% connected components.  The edge $\cs{2,4}$ is added to $\cG$ since the
% intersection of $\cs{2,5}$ and $\cs{4,5}$ is $\cs{5\right}$, a separator in
% the junction tree representation of $ \cG$.  $\cs{3,4}$ is added since it
% reduces the number of connected components.  The next edge in the
% filtration at $r=.474$ is $\cs{1,2}$, which cannot be added, since
the intersection of prime components $\cs{1,3}$ and $\cs{2,4,5}$ is empty, and
therefore not a separator.
% Finally $\cs{1,4}$ is added since $\cs{3}$ is a separator (See
% \Fig{fig:Alg1_Example}).

\begin{table}[ht]
\begin{center}
\begin{tabular}{| C C C C C C|}
\hline\Strut
\text{Coordinate}& V_1 & V_2 &V_3 &V_4 &V_5 \\
\hline\Strut
x  & 0.686 & 0.214 & 0.846 & 0.411 & 0.089 \\
y  & 0.151 & 0.194 & 0.420 & 0.567 & 0.553 \\
\hline
\multicolumn{6}{c}{}
\end{tabular}
\caption{Vertex set used to illustrate Algorithm 1.\label{tab:Alg1_Vertex}}
\end{center}
\end{table}

\begin{table}[ht]
\begin{center}
\begin{tabular}{|L C c R|}
\hline
\hfil\text{Cliques}& \text{Separators} & r & \text{Update} \\
\hline
 \cq{1}  \cq{2}\cq{3}\cq{4}\cq{5} & -  & 0 & -\quad{}  \\
 \cq{1,3}\cq{2}\cq{4}\cq{5} & - & 0.313 & \edg{1,3}\\
 \cq{1,3}\cq{2}\cq{4,5} & -  & 0.321 & \edg{4,5}  \\
 \cq{1,3}\cq{2,5}\cq{4,5} & \cq{5 }  & 0.379 & \edg{2,5}  \\
 \cq{1,3}\cq{2,4,5} & -  & 0.421 & \edg{2,4} \\
 \cq{1,3}\cq{3,4}\cq{2,4,5} & \cq{3 }\cq{4 }  & 0.459 & \edg{3,4}  \\
 \cq{1,3}\cq{3,4}\cq{2,4,5} &  \cq{3 }\cq{4 }  & 0.474 & \sim\edg{1,2} \\
 \cq{1,3,4}\cq{2,4,5} & \cq{4 }  & 0.498 & \edg{1,4} \\
\hline
%\multicolumn{4}{c}{}
\end{tabular}
\end{center}
\caption{Evolution of cliques and separators in the junction tree
  representation of $\cG$ as edges are added according to
  \Alg{filtrationalg}.  The proposed addition of edge $\edg{1,2}$ is 
  rejected.}\label{tab:Alg1_Example}
\end{table}

\subsection{Algorithm deletes few edges}

\begin{figure}[htb]
\begin{center}
  \includegraphics[height=65mm]{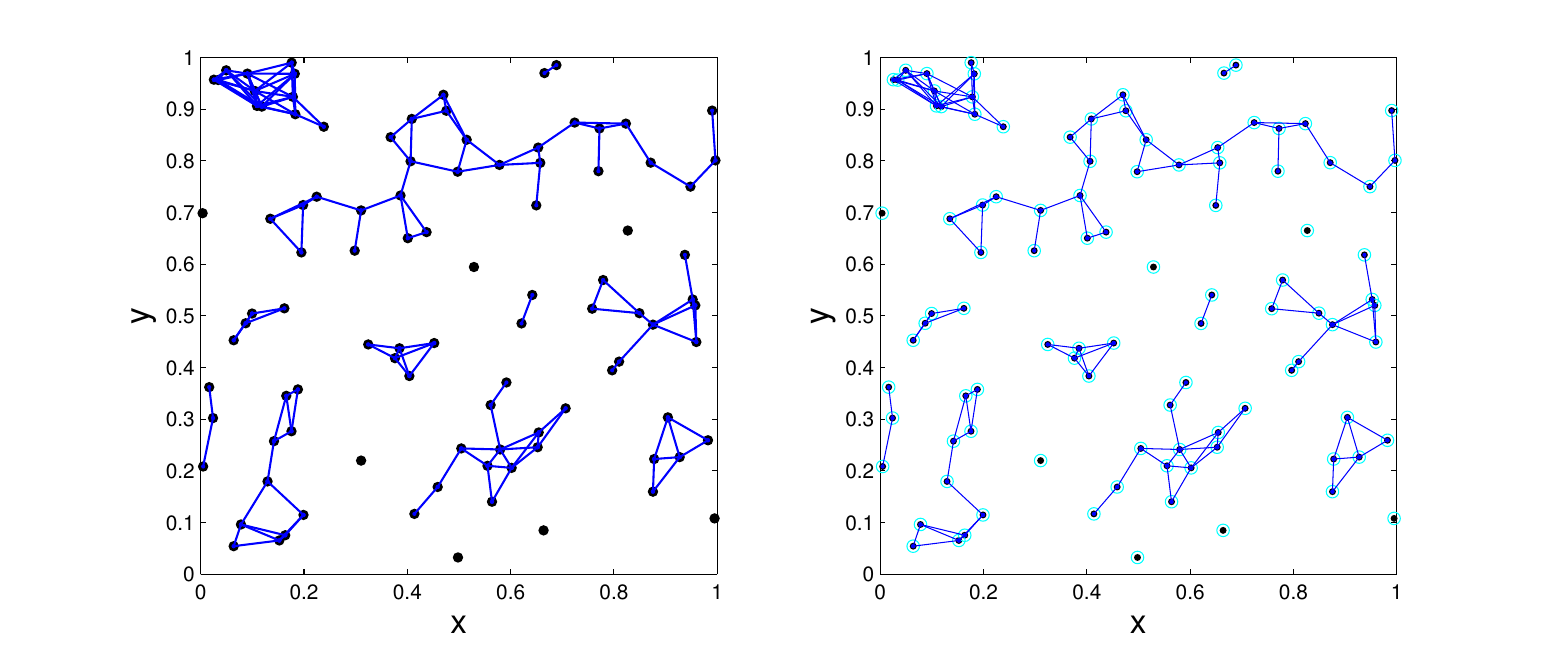}
  \caption{(a) The $1$-Skeleton of \Cech\ complex given the displayed point set and
  $r=0.05$. (b) The decomposable graph for the same complex, point set, and radius
  output by Algorithm \ref{filtrationalg}.}\label{fig:AlgExample}
\end{center}
\end{figure}

It is interesting to note that the number of proposed edge additions that are
rejected by the algorithm is typically quite small.  To illustrate this we
applied \Alg {filtrationalg} to a filtration of \Cech\ complexes
corresponding to $100$ points sampled uniformly from the unit square, with
radius $r=0.05$.  In \Fig {fig:AlgExample} the graph $\cG$ output by the
algorithm is compared to the $1$-skeleton of the \Cech\ complex (with no
decomposability restriction) with the same radius $r=0.05$.  Few edges appear
in the \Cech\ complex but not in $\cG$.  This occurs because geometric graphs
tend to be triangulated, in the sense that if edges $\edg{v_1,v_2}$ and
$\edg{v_2,v_3}$ belong to a geometric graph, then very likely the edge
$\edg{v_1,v_3}$ will also be in the graph, preserving decomposability.

 \begin{figure}[!ht]
\begin{center}
  \includegraphics[height=60mm]{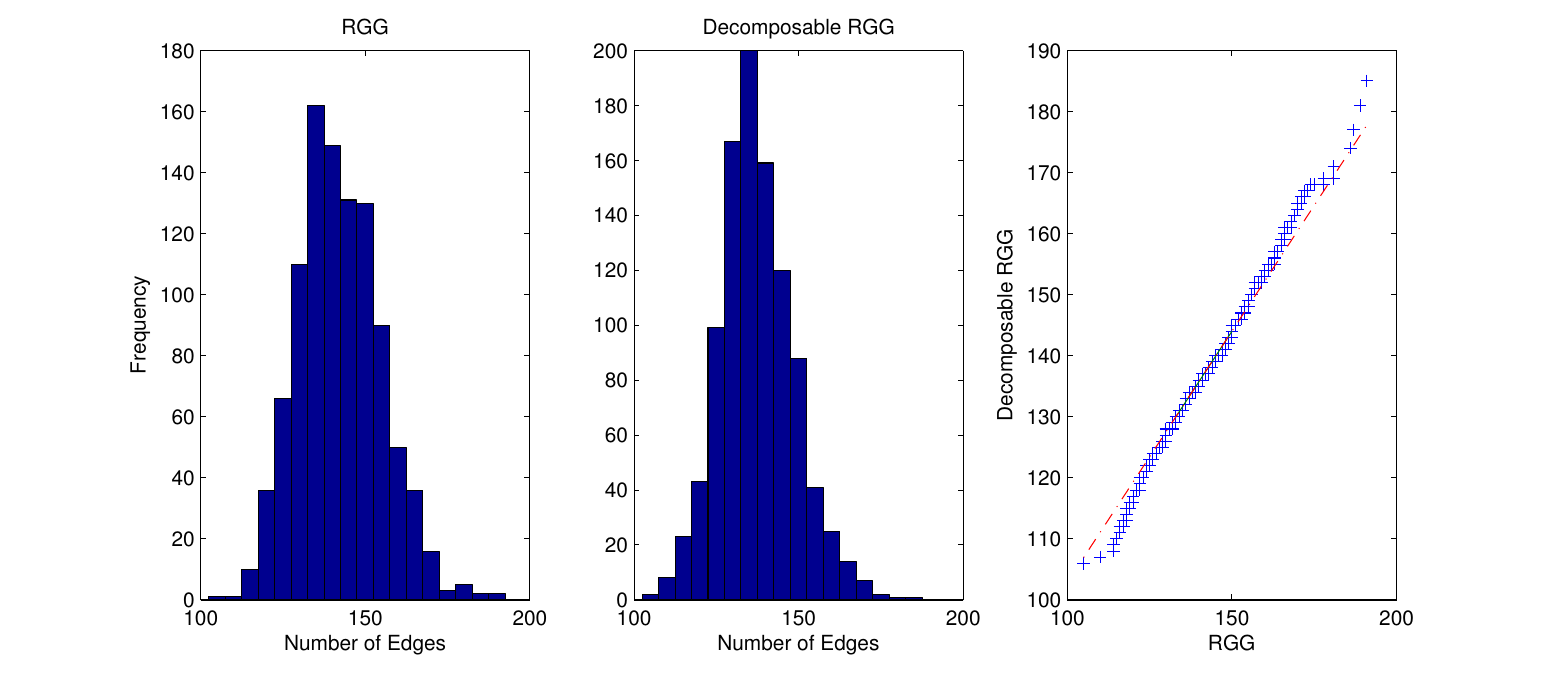}
  \caption{Distribution of edge counts for both unrestricted and decomposable
    graphs.  Graphs were computed using \Cech\ complex filtrations with
    $\Num=100$ and $V_i\iid\Un([0,1]^2)$.\label{fig:Dist_Edges}}
\end{center}
\end{figure}

\clearpage

%\cellarage
% Trim order: l b r t  (rotating -90 degrees: -> t l b r)

%\clearpage

%\clearpage
% Trim order: l b r t  (rotating -90 degrees: -> t l b r)

%\clearpage
% Trim order: l b r t  (rotating -90 degrees: -> t l b r)

%\clearpage

% Trim order: l b r t  (rotating -90 degrees: -> t l b r)
%\clearpage

%\clearpage

%\clearpage

%\clearpage

%\clearpage

% Trim order: l b r t  (rotating -90 degrees: -> t l b r)

\end{document}